\documentclass[12pt,notitlepage]{amsart}
\usepackage{latexsym,amsfonts,amssymb,amsmath,amsthm} %{pstricks,pst-node,epsf}
\DeclareFontFamily{U}{wncy}{}
\DeclareFontShape{U}{wncy}{m}{n}{<->wncyr10}{}
\DeclareSymbolFont{mcy}{U}{wncy}{m}{n}
\DeclareMathSymbol{\Sh}{\mathord}{mcy}{"58}

\newcommand{\mz}{\ensuremath{\mathbb Z}}
\newcommand{\mr}{\ensuremath{\mathbb R}}

\newcommand{\mq}{\ensuremath{\mathbb Q}}

\newcommand{\mf}{\ensuremath{\mathbb F}}

\newcommand{\mymod}{\ensuremath{\negthickspace \negmedspace \pmod}}
\newcommand{\shortmod}{\ensuremath{\negthickspace \negthickspace \negthickspace \pmod}}

\newcommand{\onehalf}{\ensuremath{ \frac{1}{2}}}
\newcommand{\half}{\ensuremath{ \frac{1}{2}}}
\newcommand{\thalf}{\textstyle \frac{1}{2}}

\newcommand{\notdivtext}{\ensuremath{\not |}}

\theoremstyle{plain}		
	\newtheorem{mytheo}{Theorem}[section]
	
	\newtheorem{myprop}[mytheo]{Proposition}
	\newtheorem{mycoro}[mytheo]{Corollary}
     \newtheorem{mylemma}[mytheo]{Lemma}
	\newtheorem{mydefi}[mytheo]{Definition}

\theoremstyle{remark}

\begin{document}

\title{On the nonvanishing of elliptic curve L-functions at the central point}
\author{Matthew P. Young} 
\address{American Institute of Mathematics, 360 Portage Ave.,
Palo Alto, CA 94306-2244}
\email{myoung@aimath.org}
\thanks{This research was partially conducted during the period the author was employed by the Clay Mathematics Institute as a Liftoff Fellow. This research was partially supported by an NSF Mathematical Sciences Post-Doctoral Fellowship.}
\begin{abstract}
We show that a large number of elliptic curve L-functions do not vanish at the central point, conditionally on the Generalized Riemann Hypothesis and on an assumption on the regular distribution of the root number.  %Some hypothesis on the root number is necessary because it has not yet been ruled out that the root number is $-1$ for almost all elliptic curves.
\end{abstract}
\maketitle

\section{Introduction and statement of results}
\subsection{Statement of main results}
There has been a lot of interest in the question of how many L-functions in a given family vanish at the central point.  There are results for many families, including Dirichlet L-functions (\cite{S}, \cite{IS1}), weight $k$ Hecke L-functions of level $N$ (\cite{IS2}, \cite{KMVdK}), and quadratic twists of a fixed elliptic curve (\cite{PP}, \cite{H-B}), to name a few examples.  In this paper we study the question of nonvanishing for the family of elliptic curves given by the Weierstrass equations
\begin{equation*}
E_{a, b}: y^2 = x^3 + ax + b
\end{equation*}
where $1 \leq a \leq X^{1/3}$, $1 \leq b \leq X^{1/2}$ (so $|\Delta| = 16(4a^3 + 27b^2) \ll X$), and $X$ is a large parameter.  Actually, we impose the additional condition that there are no primes $p$ such that $p^2 | a$ and $p^3 | b$.
This condition insures that no curve in this family is a quadratic twist of another (the Weierstrass equation of the curve $E_{a, b}$ twisted by $d$ is given by $y^2 = x^3 + ad^2x + bd^3$).  In addition, the condition insures that the above Weierstrass equation is minimal at all primes $p \neq 2$.  Let
\begin{equation*}
S = \left\{ (a,b) \in \mz^2 : (b,2) = 1 \text{ and } p^2 | a \Rightarrow p^3 \notdivtext b \right\}.
\end{equation*}
Having $b$ odd forces minimality at $p=2$ so all Weierstrass equations under consideration are global minimal models.  This fact can be deduced easily from Lemma 10.1 of \cite{Knapp} for example.  Note further that for $(a,b) \in S$, $\Delta \neq 0$, so any Weierstrass equation coming from $S$ gives an elliptic curve.

We use the standard analytic normalization of an L-function to have line of symmetry $\text{Re } s = \half$ and central point $s=\half$.

Our main result is
\begin{mytheo}
\label{thm:nonvanish}
Assume the Generalized Riemann Hypothesis and that the root number is regularly distributed as $a$ and $b$ vary in arithmetic progressions to all moduli $\ll X^{2/9 + \varepsilon}$. %\eqref{eq:Vsum}.  
Let 
\begin{equation*}
R_0(X) = \# \{(a,b) \in S: L(\thalf, E_{a,b}) \neq 0, 1 \leq a \leq X^{\frac13}, 1 \leq b \leq X^{\frac12}  \}.  
\end{equation*}
Then
\begin{equation*}
%\mathop{\sum \sum}_{\substack{
%(a, b) \in S \\ 
%X^{1/3} \leq a \leq 2X^{1/3} \\ 
%X^{1/2} \leq b \leq 2X^{1/2} \\
%L(1/2, E_{a, b}) \neq 0}} 1 
R_0(X) \gg X^{5/6 - \varepsilon}
\end{equation*}
for any $\varepsilon >0$.
\end{mytheo}

The exact requirement on the distribution of the root number has not been made precise yet for ease of exposition.  In fact, we require something which is somewhat weaker than what is stated in the Theorem but is more complicated to state and requires some additional notations.  The precise requirement is \eqref{eq:Vsum}.
%The assumption \eqref{eq:Vsum} is a technical requirement on the oscillation in sign of the root number, described in detail in the following section.

It is known, due to Kolyvagin \cite{Kolyvagin}, Gross and Zagier \cite{GZ}, and others that if the analytic rank of $L(s,E)$ is either $0$ or $1$ then the algebraic rank of $E$ equals the analytic rank and that the Tate-Shafarevich group is finite.  In particular, if $L(\thalf, E) \neq 0$ then $E$ has algebraic rank $0$.  We therefore deduce
\begin{mycoro}
\label{coro:rank0}
On the same assumptions as Theorem \ref{thm:nonvanish}, the number of rank $0$ elliptic curves with conductor $\leq X$ is $\gg X^{5/6 - \varepsilon}$ for any $\varepsilon > 0$.
\end{mycoro}

To establish this result we prove an asymptotic formula for the average behavior of a long sum of Dirichlet coefficients of elliptic curve L-functions (to be precise this sum is the `first part' of the asymmetric approximate functional equation).  See Theorem \ref{thm:LU}, which is of independent interest.

%\subsection{Statement of results}
%From the Lindel\"{o}f bound and taking an appropriate test function $w$ we deduce

\subsection{Overview of the proof}
We attack the nonvanishing problem using the standard analytic technique of comparing a lower bound for the first moment of the central values to an upper bound on the second moment of the central values.
%We now recall the standard analytic techniques for proving a nonvanishing result for a family $\mathcal{F}$ of L-functions $L(s, f)$.  
%There are essentially two ingredients, namely showing bounds of the type
Precisely, the goal is to prove bounds of the type
\begin{equation}
\tag{a}
\label{eq:a}
\sum_{f \in \mathcal{F}} L(1/2, f)\geq \mathcal{A}, 
\end{equation}
and
\begin{equation}
\tag{b}
\label{eq:b}
\sum_{f \in \mathcal{F}} L^2(1/2, f) \leq \mathcal{B}
\end{equation}
for a family $\mathcal{F}$ of L-functions.
From these two inequalities and a simple application of Cauchy's inequality we obtain
\begin{equation*}
\sum_{\substack{f \in \mathcal{F} \\ L(1/2, f) \neq 0}} 1 \geq \frac{\mathcal{A}^2}{\mathcal{B}}.
\end{equation*}
We are therefore interested in proving the bounds \eqref{eq:a} and \eqref{eq:b} for the family of all elliptic curves.

The upper bound \eqref{eq:b} with $\mathcal{B} = |\mathcal{F}|^{1 + \varepsilon}$ is easily obtained by using the Lindel\"{o}f bound; we are assuming the Generalized Riemann Hypothesis anyway for other reasons.  

The lower bound \eqref{eq:a} is the critical issue.  Establishing such a bound for the family of all elliptic curves is more subtle than for other families previously considered because the behavior of the root number is not well-understood for the family of all elliptic curves.
For the remainder of this section we discuss the problem of the lower bound.

We shall use the following expression for $L(\thalf, E)$ (the `approximate' functional equation)
\begin{equation}
\label{eq:AFE}
L(1/2, E) = \sum_{n=1}^{\infty} \frac{\lambda_E(n)}{\sqrt{n}} Y\left(\frac{2\pi n}{U} \right) 
+ \epsilon_E \sum_{n=1}^{\infty} \frac{\lambda_E(n)}{\sqrt{n}} Y\left(\frac{2\pi n}{V} \right),
\end{equation}
where $UV = N$,
\begin{equation*}
Y(u) = \frac{1}{2 \pi i} \int_{(\sigma)} u^{-t} G(t) \Gamma(1 + t) \frac{dt}{t},
\end{equation*}
$G(t)$ is an even, bounded, holomorphic function in $-4 < \text{Re }t < 4$ satisfying $G(0) = 1$, and $\epsilon_E = \pm 1$ is the root number of $E$.  Here $Y$ is a cutoff-type function because $Y(u) = 1 + O(u)$ for $u$ small and $Y(u) \ll (1 + u)^{-3}$ for $u$ large.  For technical reasons we further restrict to $G(t)$ such that $Y^{(k)}(u) \ll_{k,M} (1 + u)^{-M}$ for $u \geq 0$ and $k\geq 0$; for instance the choice $G(t) = 1$ gives $Y(u) = e^{-u}$.  As a notational shorthand we write 
\begin{equation*}
L_T = \sum_n \frac{\lambda_E(n)}{\sqrt{n}} Y\left(\frac{2\pi n}{T}\right) 
\end{equation*}
so that the approximate functional equation reads
\begin{equation*}
L(\thalf, E) = L_U + \epsilon_E L_V.  
\end{equation*}
The approximate functional equation says that we may represent the central value of $L(s, E)$ by two sums of Dirichlet coefficients, the product of whose lengths is the conductor $N$ of $E$.  In standard applications of the approximate functional equation it is usually best to take $U = V$ so as to minimize the maximum of the lengths of the two sums, but this choice is good only when the sign in the functional equation is well-understood.  In Section \ref{section:rootnumber} we elaborate upon the reasons why the root number is ill-understood for the family of all elliptic curves.  For now it suffices to say that the root number is strongly related to the variation in sign of the M\"{o}bius function of the polynomial $4a^3 + 27b^2$ as $a$ and $b$ vary.

For purposes of obtaining a lower bound of the type (a) we endeavor to minimize $V$ subject to the constraint that we can still handle the summation up to $U$ (on average over the family of course).  In this way we minimize the influence of the root number.  If we could take $U$ larger than $N$ then the root number would be effectively eliminated (more precisely, its effect would be implicitly detected by the sum $L_U$).  Unfortunately, the technology of analytic number theory does not appear to be strong enough to take $U$ this large.  

The following definition sets the parameters for our family.
\begin{mydefi}
Let %$S = \{ (a, b) \in \mz^2 : d^2 | a, d^3 | b \Rightarrow d = \pm 1 \}$.  Let 
$w \in C_0^{\infty}(\mr^{+} \times \mr^{+})$, $w \geq 0$, $A= X^{1/3}$, $B = X^{1/2}$, set
\begin{equation*}
w_X(a,b) = w\left( \frac{a}{A}, \frac{b}{B} \right),
\end{equation*}
and
\begin{equation*}
|S_X| = \sum_{(a, b) \in S} w_X(a,b).
\end{equation*}
\end{mydefi}
Thus $|S_X|$ counts the curves in our family with smooth weights.  It is simple to show
\begin{equation*}
|S_X| = \onehalf \zeta^{-1}(5) AB \widehat{w}(0, 0) \left(1 - 2^{-5} \right)^{-1} + O(B).
\end{equation*}
The following gives the average behavior of the first part of the approximate functional equation.
\begin{mytheo}  Assume the Generalized Riemann Hypothesis.
Let $\nu < 7/9$ and set $U = X^{\nu}$.  Then
%\begin{mytheo}  
%With assumptions as in Theorem \ref{thm:LUM}, we have
\label{thm:LU}
\begin{equation}
\label{eq:LUsum}
\frac{1}{|S_X|} \sum_{(a, b) \in S} L_U w_X(a,b ) \sim c_S,
\end{equation}
as $X \rightarrow \infty$, where $c_S$ is the (positive) arithmetical constant given by \eqref{eq:cS} (the quantity $Q(p^k)$ in the expression \eqref{eq:cS} is given by Definition \ref{def:Q}).
\end{mytheo}
The important feature of Theorem \ref{thm:LU} is the large allowable size for $U$; the proof becomes significantly simpler in the estimation of the remainder term if $\nu$ is reduced.  

The asymptotic \eqref{eq:LUsum} is derived by careful evaluation and estimation of a sum quite similar to
\begin{equation}
\label{eq:nsum}
\sum_{(a, b) \in S} \sum_{n \leq U} \frac{\lambda_{a, b}(n)}{\sqrt{n}} w_X(a,b).
\end{equation}
It is not obvious, yet not altogether surprising, that methods of estimation for the above sum but for $n$ ranging over primes (or prime powers) should be applicable to the above sum.  In \cite{Y} we proved a smoothed analog of the following
\begin{equation}
\label{eq:oldpsum}
\sum_a \sum_b \sum_{p \leq U}\frac{\lambda_{a, b}(p)}{\sqrt{p}} w_X(a,b) \ll (AB)^{1 - \varepsilon}
\end{equation}
for $U = X^{\nu}$, $\nu < 7/9$ and some positive $\varepsilon$ sufficiently small (with respect to $\nu$).  This result is conditional on the Generalized Riemann Hypothesis for Dirichlet L-functions.  The challenge is obtaining \eqref{eq:oldpsum} for large $U$.  

The sum \eqref{eq:nsum} is more difficult to study than \eqref{eq:oldpsum} because there is a convenient character sum expression for $\lambda_{a,b}(p)$ (see \eqref{eq:lambda} in Section \ref{section:review}) whereas $\lambda_{a,b}(n)$ is given in terms of $\lambda_{a,b}(p)$ for $p|n$ via the Hecke relations.  Obtaining results on \eqref{eq:nsum} from methods for bounding \eqref{eq:oldpsum} is technically difficult and constitutes a large portion of this paper.  It does not appear natural to prove a general result showing that a bound for \eqref{eq:oldpsum} for given $U$ leads to the analogous bound on \eqref{eq:nsum} (with the same limit on the size of $U$), yet the techniques of this paper should be applicable in general situations.  The main arguments giving this translation are in Section \ref{section:cleaning}.  The key is the symmetric-square L-function.

With the application of bounding the average rank $r$ of the family of all elliptic curves, Brumer \cite{B} obtained \eqref{eq:oldpsum} for $U \leq X^{5/9}$, giving the bound $r \leq 1/2 + 9/5 = 2.3$.  Heath-Brown \cite{H-B} showed $U = X^{2/3 - \varepsilon}$ is allowable and obtained $r \leq 1/2 + 3/2 = 2$.  The  $7/9$ result leads to $r \leq 1/2 + 9/7 = 25/14 = 1.78...$.  It is important to mention that the application of the bound \eqref{eq:oldpsum} to an upper bound for $r$ relies on the Generalized Riemann Hypothesis for all elliptic curve L-functions, so all results on the average rank discussed here are conditional on GRH.  The bound $25/14$ is interesting because it demonstrates that a positive proportion of curves have either rank $0$ or rank $1$.  

Supposing that we know that the root number is equidistributed, we still cannot conclude that a positive proportion of curves have rank $0$, since we require $r < 3/2$ to rule out $50 \%$ curves with rank $1$ and $50 \%$ rank $2$.  On the other hand, Heath-Brown showed that the average rank of a family of quadratic twists of a given elliptic curve has average rank $r \leq 3/2$ (using GRH for elliptic curve L-functions).  For this family the root number is constant when the twisting integers are held in arithmetic progressions modulo the conductor of the fixed elliptic curve.  Because of this fact Heath-Brown was able to show the average rank bound of $3/2$ holds in both subfamilies of quadratic twists where the root number is $+1$ and where the root number is $-1$.  It follows that a positive proportion of quadratic twists do not vanish at the central point; here it is necessary to have average rank less than $2$ to make this conclusion.  At present the only families known to have average rank less than $2$ are the family of all elliptic curves considered in this paper and families of quadratic twists.

The method used to deduce upper bounds on $r$ from the estimate \eqref{eq:oldpsum} is rather wasteful; it involves counting all the zeros near the central point and then simply ignoring those zeros not at the central point.  Kowalski, Michel, and VanderKam \cite{KMVdK} were able to improve the corresponding upper bound for the order of vanishing of all weight $2$ cusp forms of level $q$ (squarefree) by getting sharp lower bounds on high derivatives of central values of $L(s, f)$.  Perhaps their method could be someday applied to the family of all elliptic curves, although at this time the approach is out of reach since even the zeroth derivative is not fully understood.

Theorem \ref{thm:LU} handles the average of the first term $L_U$ in the approximate functional equation.  The second term involves the root number and is therefore inaccessible with current technology.  Setting $V_{a, b} = U^{-1} N_{a, b}$ (so $V_{a, b} \leq X^{1 - \nu}$), we make the assumption
\begin{equation}
\label{eq:Vsum}
\sum_{(a, b) \in S} \epsilon_{E_{a, b}} L_{V_{a, b}} w_X(a,b) = o(AB).
\end{equation}
Since $L_V$ is $c_S$ on average, we expect the above sum to be small based entirely on the variation of the root number.
This is the minimal assumption we require on the root number for Theorem \ref{thm:nonvanish} but it is instructive to explore other natural estimates that imply \eqref{eq:Vsum}.

By expanding the sum $L_V$ we could also make the assumption that for any $n$ we have
\begin{equation*}
\sum_{(a, b) \in S} \epsilon_{a,b} \lambda_{a, b}(n) w_X(a,b) Y \left( \frac{2\pi n}{V_{a,b}} \right) = o\left( \frac{AB}{\log^2{X}} \frac{d(n)}{\sqrt{n}} \right),
\end{equation*}
(essentially $\sqrt{n}$ savings) and deduce \eqref{eq:Vsum} on summation over $n$.  Since $\lambda_{a, b}(n)$ is periodic in $a$ and $b \pmod{n}$ (actually  it is periodic modulo the product of primes dividing $n$ because $\lambda_{a,b}(p^k)$ is a polynomial in $\lambda_{a,b}(p)$; see Section \ref{section:complete} for further discussion), we could likewise assume that the root number is evenly distributed (with power savings) as $a$ and $b$ vary over arithmetic progressions to moduli $\ll X^{2/9 + \varepsilon}$.  Since $b \asymp B = X^{1/2}$, the modulus is less than the square root of the length of $B$, so this is not an unreasonable assumption.  Of course, improvements on Theorem \ref{thm:LU} would allow us to make weaker assumption than \eqref{eq:Vsum}.

For applications to nonvanishing all we require is the bound \eqref{eq:Vsum} but we expect the true bound to be $O(AB)^{1/2 + \varepsilon}$.  

There is an additional technical issue handled implicitly by a bound such as \eqref{eq:Vsum}, namely the problem that as $a$ and $b$ vary in $S$, the conductors of the curves $E_{a,b}$ may behave somewhat irregularly.  By freezing $U$ to be a fixed power of $X$ (not depending on $a$ and $b$) we are able to freely interchange the order of summation over $a$ and $b$ and $n$ in the analysis of the first part of the approximate functional equation.
The problem of the variation of the conductor can likely be solved by using the method of \cite{Miller} (see Appendix A) of using sieve methods to take a subset of $S$ such that $D=4a^3 +27b^2$ is squarefree; in this case the conductor is essentially $D$ (up to a power of two, but this is a minor issue).  Note that taking a positive proportion of $S$ would not change the statements of Theorem \ref{thm:nonvanish} or Corollary \ref{coro:rank0}

From Theorem \ref{thm:LU} and \eqref{eq:Vsum} we immediately deduce
\begin{mytheo}
\label{thm:firstmoment}
Let notation be as in Theorem \ref{thm:LU}.  Assuming the Generalized Riemann Hypothesis and that \eqref{eq:Vsum} holds we have
\begin{equation*}
\frac{1}{|S_X|} \sum_{(a, b) \in S} L(1/2, E_{a, b}) w_X(a,b) \sim c_S
\end{equation*}
as $X \rightarrow \infty$.
\end{mytheo}

From Theorem \ref{thm:firstmoment} and the Lindel\"{o}f bound we obtain Theorem \ref{thm:nonvanish}, as described at the beginning of this section.

The application of the Lindel\"{o}f hypothesis to obtain the upper bound \eqref{eq:b} may be somewhat unsatisfying.  The problem with treating the second moment of the central values using the approximate functional equation is that it is necessary to consider sums of Fourier coefficients of length at least $N$ (note that squaring a partial sum essentially also squares the length of the partial sum).  Since we cannot handle sums this long at present (the longest sum we can handle is essentially $N^{7/9 - \varepsilon}$ as in Theorem \ref{thm:LU}) we cannot obtain an upper bound of the correct order of magnitude, at least with this method.

It would be interesting to obtain an upper bound of the type \eqref{eq:b} unconditionally.  Since only an upper bound here is desired and because of positivity of the terms it is possible to avoid the difficulties of the distribution of the root number.  This is an interesting direction for future research.  Of course any subconvexity bound $L(1/2, E) \ll N^{\frac{1}{4} - \delta}$ allows for an unconditional upper bound in \eqref{eq:b} (which, combined with Theorem \ref{thm:firstmoment}, gives $R_0(X) \gg X^{\frac{1}{3} + 2\delta}$)
but better results may be possible on average.

\subsection{Mollification}
Theorem \ref{thm:nonvanish} shows that a large number of central values do not vanish but barely fails to show that a positive proportion are nonzero.  In practice, for a general family of L-functions, it turns out be impossible to prove a positive proportion of nonvanishing using the method described in the previous section.  The reason is that $L(1/2, f)$ occasionally takes large enough values so that the best possible value of $\mathcal{B}$ is logarithmically larger than the best possible value for $\mathcal{A}$.  By introducing a mollifier $M(f)$ (an approximation to $L(1/2, f)^{-1}$), one hopes to show
\begin{equation}
\tag{a'}
\label{eq:a'}
\sum_{f \in \mathcal{F}} L(1/2, f) M(f) \gg \mathcal{A}
\end{equation}
and
\begin{equation}
\tag{b'}
\label{eq:b'}
\sum_{f \in \mathcal{F}} L^2(1/2, f) M^2(f) \ll \mathcal{A}.
\end{equation}
and prove that a positive proportion of central values are nonzero.  The existence of such a mollifier is not {\em a priori} obvious and picking a mollifier that optimizes the implied constants can be a tricky problem.  Nevertheless, based on numerous examples \cite{IS2}, \cite{KM2}, \cite{KM3} and some general conjectures \cite{CS} at the outset we take a mollifier of the form
\begin{equation}
\label{eq:mollifier}
M(E) = \sum_{m \leq M} \frac{\rho_E(m)}{\sqrt{m}} P \left(\frac{\log{M/m}}{\log{M}} \right),
\end{equation}
where $\rho_E(m)$ are determined by
\begin{equation*}
\frac{1}{L(s, E)} = \sum_{m = 1}^{\infty} \frac{\rho_E(m)}{m^s}, \qquad \text{Re }s > 1,
\end{equation*}
and $P(x)$ is a real polynomial satisfying $P(0) = 0$ and $P(1) = 1$.  The optimal mollifier for the family of all weight $2$ cusp forms of level $N$ is essentially of this form (see (2.25) of \cite{IS2}, where $P(x) = x$).

We now state the mollified analogue of Theorem \ref{thm:LU}. 
\begin{mytheo}  Assume the Generalized Riemann Hypothesis.
\label{thm:LUM}
Let $\nu < 7/9$ and set $U = X^{\nu}$.  Then for $M(E)$ given by \eqref{eq:mollifier} with $M = X^{\kappa}$, $0 < \kappa < 7/9 - \nu$, we have
\begin{equation}
\label{eq:LUMsum}
\frac{1}{|S_X|} \sum_{(a, b) \in S} L_U M(E_{a, b}) w_X(a,b) \sim \onehalf 
\end{equation}
as $X \rightarrow \infty$.
%\end{mytheo}
\end{mytheo}

The mollified analogue of \eqref{eq:Vsum} is
\begin{equation}
%\tag{\ref{eq:Vsum}'}
\label{eq:Vsum'}
\sum_{(a, b) \in S} \epsilon_{E_{a, b}} L_{V_{a, b}} M(E_{a, b}) w_X(a,b) = o(AB).
\end{equation}
Now the mollified version of Theorem \ref{thm:firstmoment} becomes
\begin{mytheo}
Let notation be as in Theorem \ref{thm:LUM}.  Assuming the Generalized Riemann Hypothesis and that \eqref{eq:Vsum'} holds we have
\begin{equation*}
\frac{1}{|S_X|} \sum_{(a, b) \in S} L(1/2, E_{a, b}) M(E_{a, b}) w_X(a,b) \sim \onehalf 
\end{equation*}
as $X \rightarrow \infty$.
\end{mytheo}

The presence of the mollifier introduces a significant technical difficulty in the proof of Theorem \ref{thm:LUM} in comparison to Theorem \ref{thm:LU}.  Nevertheless, note that Theorems \ref{thm:LU} and \ref{thm:LUM} are of the same quality in the sense that we obtain the desired asymptotic as long as the sums have length $\ll X^{7/9 - \varepsilon}$.  The barrier to obtaining a positive proportion of nonvanishing for the family of all elliptic curves arises only from the limitation on the length of the sums in these theorems.

\subsection{The second moments of partial sums of the approximate functional equation}
It is of great interest to weaken the assumption \eqref{eq:Vsum} or \eqref{eq:Vsum'}.
The central difficulty in treating the sum $\sum_{E} \varepsilon_E L_V$ appearing in \eqref{eq:Vsum} is that harmonic analysis is incompatible with the M\"{o}bius function.  It is conceivable that an eventual proof of the equidistribution of $\varepsilon_E$ (perhaps conditional on GRH, although it is not obvious how that assumption would be useful) would also show $\sum_E \varepsilon_E L_V = o(AB)$ for $V =X^\delta$ with some $\delta > 0$.  We would like to require this bound for $\delta$ as small as possible.

One obvious approach would be to prove Theorem \ref{thm:LU} with larger $U$.  Of course, any improvement in that result would lead to an improvement on the upper bound of the average rank $r$, so there is already motivation in this direction.

A natural approach to eliminating the variation of the root number is with suitable use of Cauchy's inequality.  It becomes necessary to consider sums such as $\sum_E L_V^2$ or $\sum_E L_V^2 M^2(E)$.  What follows in this section is a collection of results addressing a variety of similar sums.  Some of the results may be surprising. 

We begin with the following
\begin{myprop}
\label{prop:LVsquared}
Let $V = X^\alpha$ for $0 < \alpha < 5/18$.  Then
\begin{equation}
\label{eq:LVsquared}
\frac{1}{|S_X|} \sum_{(a,b) \in S} L_V^2 w_X(a,b) \sim c_1 \log{V}
\end{equation}
for some $c_1 > 0$ as $X \rightarrow \infty$.
\end{myprop}
The (expected) presence of an extra logarithmic factor is an obstacle.  One might attempt to eliminate this logarithmic factor by mollifying the partial sum $L_V$.  It may be surprising that the following holds.
\begin{myprop}
\label{prop:LVMsquared}
Let $V=X^{\alpha}$, $M= X^{\beta}$, with $\alpha + \beta < 5/18$, $\alpha, \beta > 0$.  Suppose $M(E)$ is of the form \eqref{eq:mollifier}.  Then
\begin{equation*}
\frac{1}{|S_X|} \sum_{(a,b) \in S} L_V^2 M^2(E_{a,b}) w_X(a,b) \sim c_2 (\log{M})^3 
\end{equation*}
for some $c_2 > 0$ (depending on $M$) as $X \rightarrow \infty$.
\end{myprop}
Rather than dampening the large values of $L_V$ the mollifier actually makes the sum substantially larger!  A partial explanation of this phenomenon is provided by
\begin{myprop}
\label{prop:Msquared}
Let $M = X^{\beta}$ for $0< \beta < 5/18$ and suppose $M(E)$ is of the form \eqref{eq:mollifier}.  Then
\begin{equation*}
\frac{1}{|S_X|} \sum_{(a,b) \in S} M^2(E_{a,b}) w_X(a,b) \sim c_3 (\log{M})^3
\end{equation*}
for some $c_3 > 0$ (depending on $M$) as $X \rightarrow \infty$.%, and where $F(x) = x \frac{d}{dx} x P(x)$.
\end{myprop}
The most natural explanation of this behavior is that the mollifier takes rather large values somewhat independently of the size of $L_V$, at least for $V$ relatively small.  Of course it should happen that when $V$ is the square root of the conductor then the asymptotic in Proposition \ref{prop:LVMsquared} should be a constant (because we expect to mollify the central values of the L-functions!).  
It is expected that a new main term should arise once sums longer than the conductor appear.

It is curious that the mollifier itself requires mollification from the central values.

Detecting the main term in Proposition \ref{prop:LVMsquared} is somewhat subtle.  The computation involves a quadruple 
integral of a ratio of products of the Riemann zeta function.  The obvious approach to compute the integral is to use Cauchy's theorem to move the lines of integration past the poles of the zeta function and compute residues.  It turns out that the main term does not arise from the residue where all the variables are zero, but instead comes from one of the lines of integration inside the critical strip.

It is not the case that a short mollifier (as in Proposition \ref{prop:Msquared}) detects a central zero of the L-function, for we have
\begin{myprop}  Let $L(s,f)$ be the L-function attached to $f \in H_2(q)$, $H_2(q)$ being a Hecke basis of weight 2 level $q$ trivial character cusp forms with $q$ prime (note each such $f$ is primitive).  Let $M(f)$ be a mollifier of the form \eqref{eq:mollifier} with $M = q^{\alpha}$, $0< \alpha < 1/2$.  Then for $\epsilon = \pm 1$ we have
\label{prop:rootnumbermollifier}
\begin{equation*}
\sum_{\substack{f \in H_2(q) \\ \epsilon_f = \epsilon}}{}^h |M(f)|^2 \sim c_M (\log{M})^3
\end{equation*}
for some $c_M > 0$ depending on $M$ only as $q \rightarrow \infty$.  The symbol $h$ indicates that each term in the summation is appropriately weighted, precisely for any $\alpha_f$
\begin{equation*}
\sum_{f \in H_2(q)}{}^h \alpha_f := \frac{\Gamma(k-1)}{(4\pi)^{k-1}} \sum_{f \in H_2(q)} \frac{\alpha_f}{||f||^2},
\end{equation*}
where $||f||$ is the Petersson norm.
\end{myprop} 
The salient feature of this result is that the same order of magnitude is obtained whether $\epsilon = 1$ or $-1$, the first case consisting (conjecturally) of L-functions that almost never vanish at the central point, and the second case consisting of L-functions that always vanish at the central point.

B. Conrey has pointed out to me that the analogue of Proposition \ref{prop:LVMsquared} for the Riemann zeta function has been discovered, but not published, by a number of researchers, including Siegel, Selberg, and Levinson.  The matter arose in attempts to prove that a positive proportion of zeros of $\zeta(s)$ lie on the critical line.
Curiously, for the case of $\zeta(s)$ the analogues of Propositions \ref{prop:LVsquared} and \ref{prop:LVMsquared} give the same order of magnitude, so in fact Proposition \ref{prop:LVMsquared} represents a new phenomenon.

As a further curiosity, we have
\begin{myprop}
\label{prop:LV1V2M1M2}
Let $V_1 = X^{\alpha_1}$, $V_2 = X^{\alpha_2}$, $M_1 = X^{\beta_1}$, $M_2 = X^{\beta_2}$, $\alpha_1 + \alpha_2 + \beta_1 + \beta_2 < 5/9$, $\alpha_1 \neq \alpha_2$, $\alpha_1, \alpha_2 \neq 0$.  Suppose $M_1(E)$ and $M_2(E)$ are of the form \eqref{eq:mollifier}.  Then
\begin{equation*}
\frac{1}{|S_X|} \sum_{(a,b) \in S} L_{V_1} L_{V_2} M_1(E_{a,b}) M_2(E_{a,b}) w_X(a,b) \sim c_4 
\end{equation*}
for some $c_4 >0$ as $X \rightarrow \infty$.  The result holds when $\beta_1$ or $\beta_2$ (or both) are zero, in which case we interpret $M(E) = 1$.
\end{myprop}

Having $\alpha_1 \neq \alpha_2$ is the crucial difference between this result and the previous ones where an extra logarithmic factor was present.

It is probable that the results in this section could be proven with larger exponents (as long as the sums have length at most $X^{7/9 - \varepsilon}$) using the same method as in the proof of Theorem \ref{thm:LUM}.  Since our intent here is to develop the main terms we prove the
results for less than maximal exponents because doing so provides a major simplification of the estimation of the remainder terms.

\subsection{Organization of the paper}
We have already shown how Theorem \ref{thm:nonvanish} follows from Theorem \ref{thm:LU}.
We prove Theorems \ref{thm:LU} and \ref{thm:LUM} in Section \ref{section:proof}; actually, we provide the details for the (harder) case of Theorem \ref{thm:LUM} and summarize the changes in Section \ref{section:changes}.

We review the necessary facts needed about elliptic curves in Section \ref{section:review}.

In Section \ref{section:rootnumber} we compute a formula for the root number of the curve $y^2 = x^3 + ax + b$ under certain restrictions on $a$ and $b$.

The proofs of Propositions \ref{prop:LVsquared}%-\ref{prop:LV1V2M1M2}
, \ref{prop:LVMsquared}, \ref{prop:Msquared}, and \ref{prop:LV1V2M1M2} 
are carried out in Section \ref{section:secondmoments}.

The proof of Proposition \ref{prop:rootnumbermollifier} is omitted; we quickly sketch the method here.  We use the Petersson formula in the form (14.65) of \cite{IK}.  The restriction $\alpha < 1/2$ leads to an analysis of the diagonal terms only (showing that the answer is independent of the sign of the functional equation).  The evaluation of the main term arising from the diagonal symbol is nearly identical to that of the proof of Proposition \ref{prop:Msquared}.

\subsection{Acknowledgements}
I would like to heartily thank Henryk Iwaniec and Brian Conrey for many long and fruitful conversations.  I also thank David Farmer for helpful comments on an earlier version of this paper and Tsz Ho Chan and Sidney Graham for some discussions on this material.  I also thank the referee for a careful reading of the paper which enhanced the clarity and precision of this work.

\section{Review of elliptic curves}
\label{section:review}
We quickly review a few relevant properties needed about elliptic curves.  The L-function $L(s, E_{a, b})$ attached to $E_{a, b}$ is given by
\begin{equation}
\label{eq:eulerproduct}
L(s, E_{a, b}) = \prod_p \left( 1 - \frac{\lambda_{a, b}(p)}{p^s} + \frac{\psi_{N}(p)}{p^{2s}} \right)^{-1},
\end{equation}
where for $p \neq 2$,
\begin{equation}
\label{eq:lambda}
\lambda_{a, b}(p) = 
- \frac{1}{\sqrt{p}} \sum_{x \shortmod{p}} \left( \frac{x^3 + ax + b}{p} \right),
\end{equation}
and $\psi_N$ is the principal Dirichlet character modulo the conductor $N$ of $E_{a, b}$ (the conductor is described in further detail below).  For $p = 2$ we have $\lambda_{a, b}(2) = 0$.  It is this concrete character sum expression for $\lambda_{a,b}(p)$ that allows one to study the behavior of the L-functions $L(s, E_{a,b})$ as $a$ and $b$ vary over `natural' sets of integers using techniques from analytic number theory (harmonic analysis, for example).

By work of Wiles, Taylor, Breuil, Conrad, and Diamond \cite{W}, \cite{TW}, \cite{BCDT}, the completed L-function
\begin{equation*}
\Lambda(s, E_{a, b}) = \left(\frac{\sqrt{N}}{2\pi} \right)^{s + \onehalf} \Gamma(s + \thalf) L(s, E_{a, b})
\end{equation*}
is entire and satisfies the functional equation
\begin{equation*}
\Lambda(1-s, E) = \epsilon_E \Lambda(s, E),
\end{equation*}
where $\epsilon_E = \pm 1$ is the root number of $E$.  The sign of the root number determines the parity of the order of vanishing of $L(1/2, E)$.  In particular, $\epsilon_E = -1$ implies $L(1/2, E) = 0$.

The order of vanishing of $L(s, E)$ has important arithmetical applications of course because of the conjecture of Birch and Swinnerton-Dyer.  %It is commonly believed that half of all elliptic curves have rank $0$ and half have rank $1$ (with rank $2$ and higher curves not constituting a positive proportion of all elliptic curves).  Thus, statistically speaking, the rank of $E$ should be controlled by the root number.

The distribution of the root number plays an important and difficult role in the nonvanishing problem.  At present there are no results that would rule out $\epsilon_{E_{a, b}} = -1$ for almost all $a, b$.  Proving a good nonvanishing result for the family of all elliptic curves therefore necessarily depends on having good control on the variation of the root number.  %In this paper we shall not attempt to prove any result on the distribution of the root number.  One of the aims of this research has been to prove a nonvanishing result subject to a minimal assumption on the root number.

The conductor is a certain divisor of the discriminant $\Delta = -16(4a^3 + 27b^2)$ that measures how bad the reduction of $E \pmod{p}$ is for each prime $p$.  There is no simple formula for the conductor but it can be efficiently computed using Tate's algorithm (\cite{Si2}, pp. 363-368).  For primes other than $2$ and $3$, the power of $p$ dividing $N$ is one if $E$ has a node $\mymod{p}$, and the power is two if $E$ has a cusp $\mymod{p}$.  We also use the terminology $E$ has additive reduction or multiplicative reduction at $p$ corresponding to whether $E$ has a cusp or a node $\mymod{p}$, respectively.  In the special case
where $(a, b) = 1$ we have
\begin{equation*}
N = 2^\alpha 3^\beta \prod_{\substack{p | \Delta \\ (p, 6) = 1}} p,
\end{equation*}
for some $2 \leq \alpha \leq 8$ and $0 \leq \beta \leq 5$, so that $N$ is essentially squarefree.  These bounds on $\alpha$ and $\beta$ are given as Theorem 10.4 of \cite{Si2}.

\section{The distribution of the root number}
\label{section:rootnumber}

The following proposition gives a concrete formula for the root number of the curve $y^2 = x^3 + ax + b$, under some light hypotheses on $a$ and $b$.

\begin{myprop}
\label{prop:rootnumber}
Let $E_{a,b} : y^2 = x^3 + ax + b$ and suppose $4a^3 + 27b^2$ is squarefree.  Then
\begin{equation}
\epsilon_{E_{a, b}} = \mu(4a^3 + 27b^2) \left(\frac{a}{3b} \right) \chi_4(b) (-1)^{a + 1} \epsilon_2,
\end{equation} 
where 
$\chi_4$ is the primitive Dirichlet character modulo $4$, $(\frac{\cdot}{\cdot})$ is the Jacobi symbol, and $\epsilon_2$ is the local root number at $2$.
\end{myprop}

Remarks. This result illustrates the difficulty of understanding the behavior of the root number. Qualitatively speaking, the distribution of the root number is essentially controlled by the M\"{o}bius function of the polynomial $4a^3 + 27b^2$.  There is currently no technology to address this problem.

\begin{proof}
Let $N_{a,b}$ be the conductor of $E_{a,b}$ and let $N_{a,b}'=4a^3 + 27b^2$ be the odd part of $N_{a,b}$.  It is known (see \cite{R} for example) that the root number is the product of the local root numbers $\epsilon_p$ at all places $p$.  At $p = \infty$ it is true that $\epsilon_{\infty} = -1$.  It is also known \cite{R} that
\begin{equation*}
\epsilon_p =
\begin{cases}
1 &  \text{if $E$ has split multiplicative reduction at $p$}, \\
- 1 & \text{if $E$ has non-split multiplicative reduction at $p$}.
\end{cases}
\end{equation*}
Recall the definition that $E$ has split multiplicative reduction at $p$ if and only if the tangent lines at the node have $\mf_p$-rational slopes.  It can be shown elementarily (see (10.10) of \cite{Knapp}) that this is equivalent to
\begin{equation*}
\epsilon_p = - \sqrt{p} \lambda_E(p).
\end{equation*}
Since $4a^3 + 27b^2$ is squarefree there is multiplicative reduction at each $p | N'$ and by the above formula we have
\begin{equation*}
\epsilon_{E_{a, b}} = - \epsilon_2 \prod_{p | 4a^3 + 27b^2} \left(- \sqrt{p} \lambda_{a, b}(p) \right) = - \epsilon_2 \sqrt{N_{a,b}'} \lambda_{a,b}(N_{a,b}') \mu(N_{a, b}').
\end{equation*}
This formula could also be derived from Atkin-Lehner theory (Theorem 9.27 of \cite{Knapp}).

Unfortunately the local root number at $2$ does not have such a simple characterization (there are many cases to consider; see \cite{Ha}).  

We can explicitly compute $\lambda(N_{a,b}')$.  Suppose $p | 4a^3 + 27b^2$.  It is clear that  there exists $\gamma$ such that $a \equiv -3\gamma^2 \pmod{p}$ and $b \equiv 2\gamma^3 \pmod{p}$ (or see Lemma \ref{lemma:parameterization}).  We easily compute
%Then we have $\left(\frac{-3a}{p} \right) = 1$, so there exists $c$ such that $a \equiv -3c^2 \pmod{p}$.  Likewise, we have $b \equiv 2c^3 \pmod{p}$ (possibly changing $c$ to $-c$ if necessary, which does not change our expression for $a$).  Then
\begin{equation}
\label{eq:degeneratecalculation}
\begin{split}
\sqrt{p} \lambda_{a, b}(p) %& = - \sum_{x \shortmod{p}} \left( \frac{x^3 + ax + b}{p} \right) \\
& = - \sum_{x \shortmod{p}} \left( \frac{x^3 -3 \gamma^2x + 2\gamma^3}{p} \right) \\
& = - \sum_{x \shortmod{p}} \left( \frac{(x-\gamma)^2(x+2\gamma)}{p} \right) \\
& = - \sum_{\substack{x \shortmod{p} \\ x \neq \gamma}} \left( \frac{x+2\gamma}{p} \right) = \left(\frac{3\gamma}{p} \right) = \left(\frac{6b}{p} \right).
\end{split}
\end{equation}

Thus, using quadratic reciprocity (and the Jacobi symbol for notation) we have
\begin{align*}
\sqrt{N_{a,b}'} \lambda_{a,b}(N_{a,b}')  = \prod_{p | 4a^3 + 27b^2} \left(\frac{6b}{p} \right) 
 = \left(\frac{6b}{4a^3 + 27b^2} \right) \\
 = \left(\frac{2}{4a^3 + 27b^2} \right)\left(\frac{3}{4a^3 + 27b^2} \right)\left(\frac{b}{4a^3 + 27b^2} \right) \\
=  \left(\frac{2}{4a^3 + 27b^2} \right)\left(\frac{4a^3 + 27b^2}{3} \right) \left(\frac{4a^3 + 27b^2}{b} \right) (-1)^{(b-1)/2} \\
 = \chi_4(b) (-1)^a \left( \frac{a}{3b} \right),
\end{align*}
which completes the proof.
\end{proof}
It is an important question in number theory to bound sums of the M\"{o}bius function of a polynomial.  The techniques of analytic number theory are not yet strong enough for the discriminant polynomial, although there has been tremendous recent progress on other polynomials that take even fewer values (such as $x^2 + y^4$ \cite{FI} and $x^3 + 2y^3$ \cite{H-B2}).

There has been some recent progress by Helfgott on the distribution of the root number for other families of elliptic curves.  For instance, he has shown that the root number of the family
\begin{equation*}
y^2 = x(x+a)(x+b)
\end{equation*}
is evenly distributed (\cite{H2}, Corollary 5.2).  His method works for families in which a certain polynomial associated with the family (essentially the conductor) is homogeneous in $a$ and $b$ and has degree $3$.  For this family the conductor is essentially $ab(a-b)$ (up to a power of two, and on certain conditions on $a$ and $b$ similar to those given in Proposition \ref{prop:rootnumber}) and the behavior of the root number is more or less equivalent to the behavior of $\mu(ab(a-b))$ as $a$ and $b$ vary.  It is interesting to investigate the nonvanishing problem for this family of elliptic curves.  One can prove the analogue of Theorem \ref{thm:LU} with $U = X^{\nu}$ for any $\nu < 2/3$.  The analogue of \eqref{eq:Vsum} requires $V$ of size $X^{1/3 + \varepsilon}$ whereas Helfgott's (unconditional) result is not sufficient to handle $V= X^{\delta}$ for any $\delta > 0$ because that would require power savings in $\sum_a \sum_b \mu(ab(a-b))$.  Of course the easier sum $\sum_{n} \mu(n)$ can be bounded with power savings using the Riemann Hypothesis (which we are freely assuming in this work) but it is unclear how to apply L-function theory to the study of $\mu(ab(a-b))$.

The distribution of the root number is well-understood for families of quadratic twists.  In this situation the root number is essentially controlled by a Dirichlet character.  On the other hand, Helfgott has shown that for a family of semistable curves the distribution of the root number is controlled by the Liouville function of the conductor \cite{H1}.  This feature makes the study of elliptic curves in families of quadratic twists significantly easier than for the family considered in this paper.

\section{The proof of Theorem \ref{thm:LUM}}
\label{section:proof}
\subsection{Outline of the proof}
The overall plan of attack on Theorem \ref{thm:LUM} is rather simple, although the details are a bit messy.  We shall write $L_U M(E)$ as a sum of Fourier coefficients as $a$, $b$, and $n$ take values over appropriate ranges ($n$ will take a variety of values with multiplicities).  It is instructive to compare our family with the family of all weight $2$ cusp forms of level $N$ (say, prime).  There one uses the Petersson formula in order to prove
\begin{equation*}
\sum_{f \in S_2(\Gamma_0(N))^*}{}^{h} \lambda_f(n) = \delta_0(n) + R_n(N),
\end{equation*}
where $R_n(N)$ is an explicit remainder term which is small, at least for $n$ and $N$ in certain ranges.  Additional saving can be obtained on averaging over $n$.  Of course we cannot use the Petersson formula for the family of all elliptic curves, but one can still show using methods of analytic number theory (completing the sum, exponential and character sum estimates, etc.) that
\begin{equation}
\label{eq:QR}
\sum_{(a, b) \in S} \frac{\lambda_{a,b}(n)}{\sqrt{n}} = \frac{Q(n)}{n^{1/2} {n^*}^2} + R(n)
\end{equation}
where
\begin{equation*}
Q(n) = \mathop{\sum_a \sum_b}_{\shortmod{n^*}} \lambda_{a, b}(n)
\end{equation*}
and $R(n)$ is a small error term (on average at least).  Here and throughout we use the convenient notation $n^* = \prod_{p | n} p$.  The difficult part of the work is showing that $\sum_{ n \leq U} R(n)$ is small for $U$ large (up to $X^{7/9 - \varepsilon}$).  Computing the main term is established using zeta function methods (expressing the sum as an integral (or multiple integrals) of a zeta function against a kernel, moving the line(s) of integration, picking up poles, etc.).  Here our methods are quite similar to Section 5 of \cite{KMVdK}.

One difficulty in establishing the formula \eqref{eq:QR} is that a usable formula for $\lambda(n)$ is only available when $n$ is squarefree.  We can get around this by writing $n = n_1 n_2$ where $(n_1, n_2) = 1$, $n_1$ squarefree and $n_2$ power-full (i.e. consisting of squares, cubes, etc.).  Then we write $\lambda(n) = \lambda(n_1) \lambda(n_2)$ and obtain significant savings in the summation over $n_1$ using the concrete formula available for squarefree $n_1$.  We are able to ensure $n_2$ is small by assuming the Riemann hypothesis for the $L(s, E)$ as well as the symmetric-square L-function attached to $L(s, E)$.  It is conceivable that this use of RH could be removed, but we are assuming RH anyway in order to take $U \leq X^{7/9 - \varepsilon}$ and the restriction on the size of $n_2$ simplifies our arguments substantially.

\subsection{Preliminary cleaning}
\label{section:cleaning}
Recall our overall goal is to asymptotically evaluate the sum
\begin{equation*}
\sum_{(a, b) \in S} L_U M(E_{a, b}) w_X(a,b).
\end{equation*}
Further recall that we consider mollifiers of the form
\begin{equation*}
M(E_{a, b}) = \sum_{m \leq M} \frac{\rho_{a,b}(m)}{\sqrt{m}} P \left( \frac{\log{M/m}}{\log{M}} \right),
\end{equation*}
where $\rho_{a,b}(m)$ is the $m$-th Dirichlet coefficient of the Dirichlet series of $1/L(s, E_{a,b})$ and $M$ is a fixed positive power of $X$.  We have the formula
\begin{equation*}
\rho_{a,b}(k) = 
\begin{cases}
\mu(m) \lambda_{a,b}(m), &\text{if $k = ml^2$, $ml$ squarefree, $(l, \Delta_{a,b}) = 1$}, \\
0, & \text{ otherwise},
\end{cases}
\end{equation*}
which can be seen by inspection of the Euler product \eqref{eq:eulerproduct} of $L(s, E)$.
Hence
\begin{equation*}
M(E_{a, b}) = \mathop{\sum \sum}_{ml^2 \leq M} \frac{\mu^2(ml) \mu(m) \lambda_{a,b}(m)}{\sqrt{m} l} \psi_{\Delta}(l) P \left( \frac{\log{M/ml^2}}{\log{M}} \right),
\end{equation*}
where $\psi_{\Delta}$ is the principal Dirichlet character modulo the absolute value of the discriminant $\Delta = -16(4a^3 + 27b^2)$.  From now on we shall often suppress the condition $ml^2 \leq M$ by setting $P(x) = 0$ for $x < 0$.  We shall use the Hecke relations
\begin{equation*}
\lambda(m) \lambda(n) = \sum_{d | (m, n)} \psi_{\Delta}(d) \lambda\left(\frac{mn}{d^2}\right).
\end{equation*}
Then we have
\begin{equation}
\label{eq:LUMexpansion}
L_U M(E) = \sum_d \sum_l \sum_m \sum_n  \frac{\lambda(mn)}{dl \sqrt{mn}}
\mu(dm) \mu^2(dlm) \psi_{\Delta}(dl) Y\left(\frac{2\pi dn}{U} \right) P \left( \frac{\log{M/dml^2}}{\log{M}} \right).
\end{equation}
Before continuing we make a useful notational
\begin{mydefi}
For an integer $n$ let $(n)_1$ denote the product of primes exactly dividing $n$ (i.e. those primes $p$ dividing $n$ such that $p^2$ does not divide $n$).  Let $(n)_2$ denote the complementary part, i.e. $(n)_2 = n/(n)_1$.
\end{mydefi}
Let $\phi: \mr^+ \rightarrow \mr$ be a smooth nonnegative function satisfying $\phi(x) =1$ for $x \leq 1/2$, $\phi(x) = 0$ for $x \geq 2$, and $\phi(x) + \phi(x^{-1}) = 1$ for all $x > 0$.  

As an initial simplification it is useful to make the following
\begin{mylemma}  
\label{lem:LUMapprox}
Set
\begin{equation*}
H_1(d,l,m,n) = Y\left(\frac{2\pi dn}{U} \right) P \left( \frac{\log{M/dml^2}}{\log{M}} \right) \phi\left(\frac{(mn)_2}{X^{\varepsilon}}\right).
\end{equation*}
On GRH we have
\begin{align*}
L_U M(E) = \sum_d \sum_l \mathop{\sum_m \sum_n}%_{(mn)_2 \leq X^{\varepsilon}}  
\frac{\lambda(mn)}{dl \sqrt{mn}}
\mu(dm) \mu^2(dlm) \psi_{\Delta}(dl) H_1(d,l,m,n) %\\
+ O(X^{-\delta}),
\end{align*}
where $\delta > 0$ depends only on $\varepsilon > 0$.
\end{mylemma}

Here the simplification is that $(mn)_2 \ll X^{\varepsilon}$.  There are a variety of ways of enforcing such a truncation; the smooth truncation using $\phi$ is particularly simple with which to work.

\begin{proof}
We begin by setting $R = X^{\varepsilon}$ and applying the identity 
\begin{equation*}
\phi\left(\frac{(mn)_2}{R}\right) + \phi\left(\frac{R}{(mn)_2}\right) =1
\end{equation*}
to \eqref{eq:LUMexpansion} to split $L_U M(E)$ into two terms, one of which is the desired main term and the other which we must show is $O(X^{-\delta})$.  Let $S$ be the sum corresponding to the term with $\phi(R/(mn)_2)$, which we shall show is $O(X^{-\delta})$.

We begin by Mellin inversion.  We shall use the following representation for $P(x) = \sum a_j x^j$ (valid for $\sigma > 0$):
\begin{equation}
P\left(\frac{\log{y}}{\log{M}}\right) = \sum_j \frac{a_j j!}{(\log{M})^j} \frac{1}{2 \pi i} \int_{(\sigma)} \frac{y^s}{s^{j + 1}} ds.
\end{equation}
This identity is easy to derive.  For $y \leq 1$ the integral vanishes by moving the line of integration arbitrarily far to the right.  For $y > 1$ we move the line to the left and pick up the pole at $s = 0$; the residue is $(\log{y})^j/j!$.  Note $a_0 = 0$ (since $P(0) = 0$) so the integrals converge absolutely.

Let $\Phi$ be the Mellin transform of $\phi$ so that
\begin{equation*}
\phi\left(\frac{R}{(mn)_2}\right) = \frac{1}{2\pi i} \int R^{-w} (mn)_2^w \Phi(w) dw.
%\phi\left(\frac{(mn)_2}{R}\right) = \frac{1}{2\pi i} \int R^w (mn_2)^{-w} \Phi(w) dw.
\end{equation*}

Using the Mellin transforms of $P$ and $\phi$ and the definition of $Y$ we have
\begin{equation}
\label{eq:LUMintegral}
S = \sum_{j} \frac{a_j j!}{(\log{M})^j} \frac{1}{(2\pi i)^3} \int \int \int \frac{M^s}{s^{j+1}} \frac{U^v}{v} R^{-w} Z(s,v,w) (2\pi)^{-v} \Gamma(1+v) G(v) \Phi(w) ds dv dw,
\end{equation}
where
\begin{equation*}
Z(s,v,w) = \sum_d \sum_l\sum_m \sum_n  \frac{\lambda(mn) \mu(dm) \mu^2(dlm) \psi_{\Delta}(dl)}{d^{1 + s + v} l^{1 + 2s} m^{1/2 + s} n^{1/2 + v} (mn)_2^{-w}}.
\end{equation*}
In multiplicative notation we have
\begin{equation*}
Z(s, v, w) = \prod_p \mathop{\sum_d \sum_l \sum_m \sum_n }_{d + l + m \leq 1} \frac{\lambda(p^{m+n}) \psi_{\Delta}(p^{d+l}) (-1)^{d+m} \delta_2(p^{m+n}) }{p^{ d(1 + s + v) + l(1+2s) + m(1/2 + s) + n(1/2 + v)}},
\end{equation*}
where
\begin{equation*} 
\delta_2(p^{m+n}) =
\begin{cases}
1  \quad & \text{if } m+n \leq 1, \\
p^{w(m+n)} \quad & \text{if } m+n \geq 2.
\end{cases}
\end{equation*}
%Setting 
%\begin{equation*}
%L^*(s, E) = \prod_p \left(1 + \frac{\lambda(p)}{p^s} \right)
%\end{equation*}
%and pulling out the lower degree powers of $p$, we have
We need to analyze the behavior of $Z$ when the variables are near $0$, so we extract the lower-degree factors to obtain
\begin{equation*}
Z(s, v, w) = \prod_p \left(1 - \frac{\psi_\Delta(p)}{p^{1 + s + v}} + \frac{\psi_\Delta(p)}{p^{1 + 2s}} - \frac{\lambda(p)}{p^{1/2 + s}} + \frac{\lambda(p)}{p^{1/2 + v}} - \frac{\lambda(p^2)}{p^{1 + s + v -2w}} + \frac{\lambda(p^2)}{p^{1 + 2 v -2w}}\right) \eta(s,v,w),
%Z(s,v) = \frac{L(1/2 + v, E) \zeta(1 + 2v) L(1 + 2v, sym^2 E)}{L(1/2 + s, E)\zeta(1 + s + v) L(1 + s+v, sym^2 E)} \eta(s, v)
%Z(s, v) = \frac{L^*(1/2 + v, E)}{L^*(1/2 + s, E)} \frac{\zeta(1 + 2v)}{\zeta(1 + s + v)} \frac{L(1 + 2v, sym^2 E)}{L(1 + s + v, sym^2 E)} \eta(s, v),
\end{equation*}
where $\eta(s,v,w)$ is given by an Euler product that is absolutely convergent for $\text{Re } s > -1/6$, $\text{Re } v > -1/6$, $\text{Re } (s - w) > -1/6$, and $\text{Re } (v - w) > -1/6$.  In what follows we write $\eta_i$, $i=1,2,3$, for an Euler product satisfying the same convergence conditions as $\eta$.
%Here $L(s, sym^2 E)$ is the symmetric-square L-function associated to $L(s, E)$; the coefficient of $p^{-s}$ in the Dirichlet series expansion of $L(s, sym^2 E)$ is $\lambda_E(p^2)$.

%We shall use the following estimates (valid for $\text{Re } s > 0$ and $X \gg N$)
%\begin{equation*}
%L(1+s, sym^2 E)^{\pm 1} = \left[\prod_{p \leq X^{\varepsilon}} \left(1 + \frac{\lambda_E(p^2)}{p^{1+s}} \right)^{\pm 1} \right] \eta_{\pm}(s) %\left(1 + O(X^{-\delta}) \right),
%\end{equation*}
%%\begin{equation*}
%%L(1+s, sym^2 E)^{\pm 1} = \sum_{n \ll X^{\varepsilon}} \frac{\lambda_{sym^2 E}(n)}{n^{1 + s}} + O(X^{-\delta}),
%%\end{equation*}
%uniformly with respect to $E$ with conductor $N \ll X$.
%These are easily deduced by the Riemann hypothesis for $L(s, sym^2 E)$, although clearly we could assume something weaker.  The point is that %we may represent the symmetric-square L-function by a uniformly short sum at the boundary region of absolute convergence of the Euler% %product.  We may likewise restrict the primes present in the Euler product of $\eta_\pm(s)$ and in $\eta(s, v)$.
%Define $\gamma(s,v,w)$ by $Z(s,v,w) = \gamma(s,v,w) \eta(s,v,w)$.  
By factorizing (and using $\lambda(p^2) = \lambda^2(p) - \psi_{\Delta}(p)$ to simplify) we have
\begin{align*}
Z = &\prod_p \left(1 - \frac{\lambda(p)}{p^{1/2 + s}} + \frac{\psi_\Delta(p)}{p^{1 + 2s}} \right)
\left(1 + \frac{\lambda(p)}{p^{1/2 + v}} + \frac{\lambda(p^2)}{p^{1 + s + v}} - \frac{\lambda(p^2)}{p^{1 + s + v -2w}} + \frac{\lambda(p^2)}{p^{1 + 2 v -2w}} \right) \eta_1 
\\
 = &\prod_p \left(1 - \frac{\lambda(p)}{p^{1/2 + s}} + \frac{\psi_\Delta(p)}{p^{1 + 2s}} \right) \left(1 + \frac{\lambda(p)}{p^{1/2 + v}} + \frac{\lambda(p^2)}{p^{1 + 2v}} \right) 
\left(1 - \frac{\lambda(p^2)}{p^{1 + 2v}} \right) 
\left(1 + \frac{\lambda(p^2)}{p^{1 + s+v}} \right) \\
&  \cdot \prod_p \left(1 - \frac{\lambda(p^2)}{p^{1 + s+ v -2w}} \right) 
\left(1 + \frac{\lambda(p^2)}{p^{1 + 2v -2w}} \right) \eta_2 \\
 = & \frac{L(1/2 + v, E) L(1 + s + v, sym^2 E) L(1 + 2v - 2w, sym^2 E)}{L(1/2 + s, E) L(1 + s +v -2w, sym^2 E) L(1 + 2v, sym^2 E)} \eta_3.
\end{align*}

%Let
%\begin{equation*}
%\gamma(s, v) = \frac{\zeta(1 + 2v)}{\zeta(1 + s + v)} \frac{L(1 + 2v, sym^2 E)}{L(1 + s + v, sym^2 E)} \eta(s, v).
%\end{equation*} 
%We now simply truncate the Dirichlet series expansion of $\gamma$ so that the only integers appearing are $\ll X^{\varepsilon}$.  Let $\gamma_t$ be this truncated Dirichlet series.  By the Riemann hypothesis for the symmetric-square L-functions %and the bound $L(1, sym^2E)^{-1} \ll %X^{\varepsilon}$ we have
%\begin{equation*}
%\gamma(s, v) = \gamma_t(s,v) + O(X^{-\delta}),
%\end{equation*}
%valid for $\text{Re } s > 0$, $\text{Re } v > 0$.  Here $\delta$ depends on $\varepsilon$ only and the implied constant depends on $\text{Re } s$ and $\text{Re } v$.

Now we move the lines of integration in \eqref{eq:LUMintegral} to $\text{Re } s = \sigma$, $\text{Re } v = \sigma$, and $\text{Re } w = 1/6$, where $\sigma > 0$ is small.  Then we use the Riemann Hypothesis to get the bounds $L^{\pm 1} (1/2 + \sigma + it) \ll_{\epsilon'} (|t| X)^{\epsilon'}$ and $L^{\pm 1}(2/3 + 2 \sigma + it) \ll (|t| X)^{\varepsilon'}$ where $L(s)$ is either $L(s, E)$ or $L(s, sym^2 E)$.  Here $\epsilon' > 0$ is arbitrarily small and the implied constants depend only on $\epsilon'$.  Inserting these bounds into the integral we obtain
\begin{equation*}
S \ll X^{\epsilon' + \sigma - \epsilon/6}.
\end{equation*}
Taking $\sigma$ and $\epsilon'$ sufficiently small proves the result with any $\delta < \varepsilon/6$.
%The conclusion is that in the expression \eqref{eq:LUMintegral} we may replace $Z(s, v)$ by a modified zeta-function, say $Z_\varepsilon(s, %v)$ with an error term $O(X^{-\delta})$.  Here $Z_\varepsilon$ is given by the same formula as $Z$ except with $\lambda$ replaced with the %multiplicative function $\lambda_\varepsilon$, where 
%\begin{equation*}
%\lambda_{\varepsilon}(p^k) =
%\begin{cases}
%\lambda(p),& \qquad k=1 \\
%\lambda(p^k),& \qquad p^k \leq X^{\varepsilon}, k \geq 2 \\
%0,& \qquad \text{otherwise}.
%\end{cases}
%\end{equation*}
%Now we simply reverse the steps that led from the sum to the integral representation of $L_U M(E)$.  Replacing $\gamma$ by $\gamma_t$ is equivalent to the restriction $(mn)_2 \ll X^{\varepsilon}$ (up to $O(X^{-\delta})$ error), which compeletes the proof.
\end{proof}

As an additional simplification we have
\begin{mylemma}
\label{lem:LUMavgclean}
On GRH we have
\begin{multline}
\label{eq:LUMapprox}
\sum_{(a, b) \in S} L_U M(E_{a, b}) w_X(a,b) \\
= \sum_{d} \sum_{l} \sum_{\substack{c | dl \\ c \leq X^{\varepsilon}}} \mu(c) \mathop{\sum_{m} \sum_{n}}_{\substack{(m, dl) = 1}} \mathop{\sum \sum}_{\substack{(a,b) \in S \\ \Delta \equiv 0 \shortmod{c}}} \frac{\lambda_{a,b}(mn) \mu^2(dlm) \mu(dm)}{dl\sqrt{mn}} H(a,b,d,l,m,n) \\
+ O(AB X^{-\delta}),
\end{multline}
for some $\delta > 0$ (depending on $\varepsilon > 0$), where
\begin{equation*}
H(a,b,d,l,m,n) = Y\left(\frac{2\pi dn}{U} \right) P \left( \frac{\log{M/dml^2}}{\log{M}} \right) \phi\left(\frac{(mn)_2}{X^{\varepsilon}}\right) w_X(a,b).
\end{equation*}
\end{mylemma}

Here the simplification is that the character $\psi_{\Delta}(dl)$ is eliminated via M\"{o}bius inversion and we have restricted the new modulus $c$ to be small.
 
\begin{proof}
Set
\begin{equation*}
T_{dl}(a, b) = \mathop{\sum_{m} \sum_{n}}_{\substack{%(mn)_2 \ll X^{\varepsilon} \\ 
(m, dl) = 1}}
\frac{\lambda_{a,b}(mn) \mu(m)}{\sqrt{mn}} Y\left(\frac{2\pi dn}{U} \right) P \left( \frac{\log{M/dml^2}}{\log{M}} \right) \phi\left(\frac{(mn)_2}{X^{\varepsilon}}\right),
\end{equation*}
and
\begin{equation*}
T(a, b) = \sum_{d} \sum_l \frac{\psi_{\Delta}(dl) \mu^2(dl) \mu(d)}{dl} T_{dl}(a,b).
\end{equation*}
Then Lemma \ref{lem:LUMapprox} reads
\begin{equation*}
L_U M(E_{a,b}) = T(a, b) + O(X^{-\delta}).
\end{equation*}

By the Riemann hypothesis for $L(s, E_{a,b})$ we see that (using the contour integration methods in the proof of Lemma \ref{lem:LUMapprox}, for instance)
\begin{equation}
\label{eq:Tdlbound}
T_{d l}(a, b) \ll X^{\varepsilon'},
\end{equation}
for any $\varepsilon' > 0$ (we shall take $\varepsilon' = \varepsilon/2$ below) uniformly for $dl^2 \leq M, a \ll A$, and $b \ll B$.

By M\"{o}bius inversion we write
\begin{equation*}
\psi_{\Delta}(dl) = \sum_{\substack{c | (\Delta, dl) \\ c \leq X^{\varepsilon}}} \mu(c) + \sum_{\substack{c | (\Delta, dl) \\ c > X^{\varepsilon}}} \mu(c).
\end{equation*}
Using this formula for $\psi_\Delta$ we write $T(a,b) = T_1 + T_2$ in the obvious way.  By rearranging the order of summation we easily have
\begin{equation*}
\mathop{\sum \sum}_{(a,b) \in S} T_1(a,b) w_X(a,b)
\end{equation*}
giving the main term in Lemma \ref{lem:LUMavgclean}.

We presently show that $T_2$ is small on average.  We compute
\begin{align*}
\mathop{\sum \sum}_{(a,b) \in S} T_2(a,b) w_X(a,b) 
= 
\mathop{\sum_d \sum_l}_{dl^2 \leq M} \frac{\mu^2(dl) \mu(d)}{dl} \sum_{\substack{c | dl \\ c > X^{\varepsilon}}} 
\mathop{\sum \sum}_{\substack{(a,b) \in S \\ \Delta \equiv 0 \shortmod{c}}} T_{dl}(a,b) w_X(a,b).
\end{align*}
Here we save in the summation over $a$ and $b$.  For any given $a$ we have $b$ restricted by $27b^2 \equiv -4a^3 \pmod{c}$.  There are $\ll d(c)$ solutions to this equation modulo $c$.  Hence we have the bound (using $\varepsilon' = \varepsilon/2$ in \eqref{eq:Tdlbound})
\begin{align*}
\mathop{\sum \sum}_{(a,b) \in S} T_2(a,b) w_X(a,b)
& \ll A X^{\varepsilon/2} \mathop{\sum_d \sum_l}_{dl^2 \leq M} \frac{1}{dl}\sum_{\substack{c | dl \\ c > X^{\varepsilon}}} d(c) \left(1 + \frac{B}{c} \right) \\
& \ll AB X^{-\varepsilon/4},
\end{align*}
which completes the proof.
\end{proof}

\subsection{Complete sum calculations}
\label{section:complete}
In order to continue our analysis it is necessary to consider sums of the form $\sum_a \sum_b \lambda_{a,b}(n)$ for $n \in \mathbb{N}$ and $a$ and $b$ ranging over certain intervals.  It is an important point that $\lambda_{a,b}(n)$ is periodic in $a$ and $b$ modulo $n^*$ (recall $n^*$ is the product of primes dividing $n$).  To elaborate, it is easily deduced from \eqref{eq:lambda} that $\lambda_{a,b}(p)$ is periodic in $a$ and $b$ modulo $p$ (at least for $p \neq 2$; periodicity of $\lambda_{a,b}(2)$ modulo $2$ is trivial since $\lambda_{a,b}(2) = 0$ for all $(a,b) \in S$).  Periodicity modulo $n^*$ follows from multiplicativity and the fact that $\lambda_{a,b}(p^k)$ is a polynomial in $\lambda_{a,b}(p)$ with coefficients periodic in $a$ and $b$ modulo $p$ (the polynomial is explicitly given by a Tchebyshev polynomial when $(p, \Delta) = 1$ whereas if $p | \Delta$ then $\lambda_{a,b}(p^k) = \lambda_{a,b}(p)^k$).

By using Poisson summation in $a$ and $b$ modulo $n^*$ we are led naturally to summing $\lambda_{\alpha, \beta}(n)$ twisted by an arbitrary additive character as $\alpha$ and $\beta$ vary over a complete set of residue classes modulo $n^*$.  Actually we require some variations where $\Delta$ has some divisibility properties.  In this section we explicitly compute these complete character sums that will arise during our averaging.  

We begin by recording the following result due to Gauss.
\begin{mylemma}
\label{lem:Gauss}
Let $r$ be odd and squarefree.  Set
\begin{equation*}
G_k(r) = \sum_{y \shortmod{r}} \left(\frac{y}{r} \right) e\left(\frac{ky}{r} \right).
\end{equation*}
Then for any $k \in \mz$,
\begin{equation*}
G_k(r) = \varepsilon_r \sqrt{r} \left(\frac{k}{r} \right),
\end{equation*}
where
\begin{equation*}
\varepsilon_r = 
\begin{cases}
1, \qquad &r \equiv 1 \pmod{4}, \\
i, \qquad &r \equiv 3 \pmod{4}.
\end{cases}
\end{equation*}
\end{mylemma}
For a proof we refer to \cite{IK}, Theorem 3.3.

Next, we have
\begin{mylemma}
\label{lem:maincharsum}
Let $r$ be odd and squarefree and let $h$ and $k$ be integers.  Then
\begin{equation*}
\sum_{\alpha \shortmod{r}} \sum_{\beta \shortmod{r}} \lambda_{\alpha, \beta}(r) e \left(\frac{\alpha h + \beta k}{r} \right)
=
\varepsilon_r \mu(r) r \left(\frac{k}{r} \right) e\left( \frac{-h^3 \overline{k}^2}{r} \right),
\end{equation*}
where $k \overline{k} \equiv 1 \pmod{r}$ provided $(k, r) = 1$.  The formula holds in case $(k, r) \neq 1$, in which case we set $\overline{k} = 0$; in this situation both sides of the identity vanish.
\end{mylemma}
Note that when $r=1$ both sides of the identity equal $1$ (including the case $k=0$ since $\left(\frac{0}{1}\right) = 1$).

In general in what follows the meaning of $\overline{n}$ should be clear from context, that is we do not necessarily always explicitly state the modulus $q$ such that $n \overline{n} \equiv 1 \pmod{q}$.  Furthermore, in case $(n, q) \neq 1$ we set $\overline{n} = 0$, but this choice is arbitrary because in practice whenever this situation occurs we also have a term of the form $\left( \frac{n}{q} \right)$, which of course vanishes.
\begin{proof}
Let $T=T(h,k, r)$ be the sum to be computed.  To begin, we claim
\begin{equation*}
\lambda_{\alpha, \beta}(r) = \frac{\mu(r)}{\sqrt{r}} \sum_{x \shortmod{r}} \left(\frac{x^3 + \alpha x + \beta}{r} \right).
\end{equation*}
This equality is easily deduced by the Chinese Remainder Theorem and induction on the number of primes dividing $r$.

Inserting this expression for $\lambda_{\alpha, \beta}(r)$ into $T$ and applying the change of variables $\beta \rightarrow \beta - \alpha x - x^3$ we obtain
\begin{align*}
T 
& =  
\frac{\mu(r)}{\sqrt{r}} \sum_x \sum_\alpha  \sum_\beta \left(\frac{\beta}{r}\right) e \left(\frac{\beta k}{r} \right)
e \left(\frac{\alpha(h - x k)}{r} \right) e \left(\frac{-x^3 k}{r} \right) \\
& =  \varepsilon_r \mu(r) r \left(\frac{k}{r} \right) \sum_{x: h \equiv xk \shortmod{r}} e \left(\frac{-x^3 k}{r} \right)\\
& = \varepsilon_r \mu(r) r \left(\frac{k}{r} \right) e\left( \frac{-h^3 \overline{k}^2}{r} \right).
\qedhere
\end{align*}
\end{proof}

It will be necessary to compute some `degenerate' sums, i.e. sums where $4\alpha^3 + 27 \beta^2 \equiv 0 \pmod{p}$.  To aid in our calculations we state the following
\begin{mylemma}
\label{lemma:parameterization}
Let $r$ be odd and squarefree.  Then every solution $(\alpha, \beta) \pmod{r}$ to $$4\alpha^3 + 27 \beta^2 \equiv 0 \pmod{r}$$ is of the form $(-3\gamma^2, 2\gamma^3)$ as $\gamma$ varies $\pmod{r}$.
\end{mylemma}
\begin{proof}
First consider the case of prime $r$.  Then it is easy to see that every solution is of the stated form.  Uniqueness also follows easily because if $\gamma^3 \equiv \gamma'^3$ and $\gamma^2 \equiv \gamma'^2 \pmod{p}$ then $\gamma \equiv \gamma' \pmod{p}$.  The case $p = 3$ requires a separate argument, but that case is trivial.  The general case for squarefree $r$ follows from the Chinese Remainder Theorem.
%: We write
%\begin{equation*}
%\alpha = \sum_{p | r} \alpha_p (r/p) (\overline{r/p}),
%\end{equation*}
%where $\alpha_p$ varies modulo $p$ and $(r/p) (\overline{r/p}) \equiv 1 \pmod{p}$.  We similarly separate $\beta$.  By the above arguments we have $\alpha_p \equiv -3 \gamma_p^2 \pmod{p}$ and $\beta_p \equiv 2 \gamma_p^3 \pmod{p}$.  Of course we may join the $\gamma_p$'s via $\gamma \equiv \sum_{p | r} \gamma_p (r/p) (\overline{r/p}) \pmod{r}$, whence $\gamma$ takes each value $\pmod{r}$ exactly once.  Then
%\begin{align*}
%\alpha & \equiv \sum_{p | r} -3 \gamma_p^2 (r/p) (\overline{r/p}) \pmod{r} \\
%& \equiv -3 \sum_{p | r} (\gamma_p (r/p) (\overline{r/p}))^2 \pmod{r} \\
%& \equiv -3 \left(\sum_{p | r} \gamma_p (r/p) (\overline{r/p})\right)^2 \pmod{r} \\
%& \equiv -3 \gamma^2 \pmod{r}.
%\end{align*}
%Of course the same argument gives $\beta \equiv 2 \gamma^3 \pmod{r}$.
\end{proof}
We begin with the simplest case (of a degenerate sum) which follows immediately from the preceding Lemma.
\begin{mycoro}
\label{lem:degenerateexpsum}
Let $r$ be odd and squarefree and set $D = 4\alpha^3 + 27\beta^2$.  Then
\begin{equation*}
\mathop{\sum_{\alpha \shortmod{r}} \sum_{\beta \shortmod{r}} }_{D \equiv 0 \shortmod{r}} 
e \left(\frac{\alpha h  + \beta k }{r} \right)
= 
\sum_{\gamma \shortmod{r}} e \left(\frac{-3\gamma^2 h  + 2\gamma^3 k }{r} \right).
\end{equation*}
\end{mycoro}
%\begin{proof}
%Using the Chinese Remainder Theorem we factor the sum as follows
%\begin{equation*}
%\mathop{\sum_{\alpha \shortmod{r}} \sum_{\beta \shortmod{r}} }_{D \equiv 0 \shortmod{r}} 
%e \left(\frac{\alpha h  + \beta k }{r} \right)
%= \prod_{p | r} \mathop{\sum_{\alpha \shortmod{p}} \sum_{\beta \shortmod{p}} }_{D \equiv 0 \shortmod{p}} 
%e \left(\frac{\alpha h (\overline{r/p}) + \beta k (\overline{r/p})}{p} \right).
%\end{equation*}
%The result immediately follows from the fact that the solutions to $4 \alpha^3 + 27\beta^2 \equiv 0 \pmod{r}$ are %parameterized by $\alpha \equiv -3\gamma^2 \pmod{r}$, $\beta \equiv 2 \gamma^3 \pmod{r}$, which we prove presently.  
%\end{proof}

To cover the case where $\lambda_{\alpha, \beta}(r)$ appears, we have
\begin{mylemma}
\label{lem:degeneratecharsum}
Let $r$ be odd and squarefree and set $D = 4\alpha^3 + 27\beta^2$.  %Suppose $(u, r) = (v, r) =1$.  
Then
\begin{equation*}
\mathop{\sum_{\alpha \shortmod{r}} \sum_{\beta \shortmod{r}}}_{D \equiv 0 \shortmod{r}} \lambda_{\alpha, \beta}(r) e \left(\frac{\alpha h  + \beta k }{r} \right)
= \frac{1}{\sqrt{r}} \left(\frac{3}{r} \right) \sum_{\gamma \shortmod{r}} \left(\frac{\gamma}{r} \right)
e\left( \frac{-3\gamma^2h + 2\gamma^3 k}{r} \right).
\end{equation*}
\end{mylemma}
In particular, the sum is zero if $3 | r$.
\begin{proof}
We use Lemma \ref{lemma:parameterization} to parameterize the solutions of $D \equiv 0 \pmod{r}$.  We can explicitly compute that
\begin{equation*}
\lambda_{-3 \gamma^2, 2\gamma^3}(r) = \frac{1}{\sqrt{r}} \left(\frac{3\gamma}{r}\right),
\end{equation*}
using the same calculation as \eqref{eq:degeneratecalculation}.
%We conclude the proof by using the multiplicativity of $\lambda(r)$.
\end{proof}

We have also
\begin{mycoro}
\label{cor:maincompletesum}
Let $r$ be odd and squarefree, suppose $t | r$, and set $r = r_0 t$. Then
\begin{multline*}
\mathop{\sum_{\alpha \shortmod{r}} \sum_{\beta \shortmod{r}}}_{D \equiv 0 \shortmod{t}} \lambda_{\alpha, \beta}(r) e \left(\frac{\alpha h  + \beta k }{r} \right) \\
= 
\varepsilon_{r_0} \mu(r_0) \frac{r_0}{\sqrt{t}} \left(\frac{3}{t} \right) \left(\frac{kt}{r_0} \right)
e\left( \frac{-h^3 \overline{k}^2 \overline{t}}{r_0} \right)
\sum_{\gamma \shortmod{t}} \left(\frac{\gamma}{t} \right)
e\left( \frac{(-3\gamma^2h + 2\gamma^3 k)\overline{r_0}}{t} \right).
\end{multline*}
\end{mycoro}
\begin{proof}
Using the Chinese Remainder Theorem we write
\begin{equation*}
\alpha = \alpha_0 t \overline{t} + \alpha_1 r_0 \overline{r_0},
\end{equation*}
where $\alpha_0$ takes values modulo $r_0$ and $\alpha_1$ takes values modulo $t$, and $t\overline{t} \equiv 1 \pmod{r_0}$ and $r_0 \overline{r_0} \equiv 1 \pmod{t}$.  We similarly separate the moduli with $\beta$.  Then our sum separates into the following product
\begin{equation*}
\sum_{\alpha_0 \shortmod{r_0}} \sum_{\beta_0 \shortmod{r_0}} \lambda_{\alpha_0, \beta_0}(r_0)  \left(\frac{(\alpha_0 h  + \beta_0 k )\overline{t} }{r_0} \right)
\mathop{\sum_{\alpha_1 \shortmod{t}} \sum_{\beta_1 \shortmod{t}}}_{D \equiv 0 \shortmod{t}} 
\lambda_{\alpha_1, \beta_1}(t) \left(\frac{(\alpha_1 h  + \beta_1 k )\overline{r_0} }{t} \right).
\end{equation*}
Applying Lemmas \ref{lem:maincharsum} and \ref{lem:degeneratecharsum} completes the proof.
\end{proof}

It is useful to make the following
\begin{mydefi}
\label{def:Q}
For any $r$ and $t$ set
\begin{equation*}
Q_{t}(r) = \mathop{\sum_{\alpha \shortmod{r^*}} \sum_{\beta \shortmod{r^*}}}_{\Delta \equiv 0 \shortmod{(r^*,t)}} \lambda_{\alpha, \beta}(r).
\end{equation*}
For $t = 1$ we set $Q(r) = Q_1(r)$.
\end{mydefi}
\begin{mylemma}
\label{lem:Q(p^2)}
For $r$ squarefree we have
\begin{equation*}
Q(r^2) = 0.
\end{equation*}
\end{mylemma}
\begin{proof}
It suffices to show $Q(p^2) = 0$.  This is trivial for $p=2$ since $E_{a,b}$ has additive reduction at $p = 2$ for any $a$,$b \in \mz$, so assume $p \neq 2$ in what follows.  By the Hecke relations,
\begin{equation*}
\lambda_{\alpha, \beta}(p^2) = \lambda_{\alpha, \beta}^2(p) - \psi_{\Delta}(p).
\end{equation*}
Hence
\begin{equation*}
Q(p^2) = \sum_{\alpha \shortmod{p}} \sum_{\beta \shortmod{p}} \lambda_{\alpha, \beta}^2(p) - p(p-1),
\end{equation*}
because there are exactly $p$ pairs $(\alpha, \beta)$ such that $\psi_{\Delta}(p) = 0$ (which follows from Lemma \ref{lemma:parameterization}).  We compute
\begin{align*}
\mathop{\sum \sum}_{\alpha, \beta \shortmod{p}}  \lambda_{\alpha, \beta}^2(p) & =
\frac{1}{p} \sum_{\alpha \shortmod{p}} \sum_{\beta \shortmod{p}} \sum_{x \shortmod{p}} \sum_{y \shortmod{p}} \left(\frac{x^3 + \alpha x + \beta}{p} \right) \left(\frac{y^3 + \alpha y + \beta}{p} \right) \\
& = \frac{1}{p} \sum_{\alpha \shortmod{p}} \sum_{\beta \shortmod{p}} \sum_{x \shortmod{p}} \sum_{y \shortmod{p}}
\left(\frac{\beta}{p} \right) \left(\frac{(x^3 - y^3) + \alpha (x- y) + \beta}{p} \right),
\end{align*}
by the change of variables $\beta \rightarrow \beta - \alpha y - y^3$.  The summation over $\alpha$ is zero unless $x = y$.  Hence this sum is $p(p-1)$, which completes the proof.
\end{proof}

\begin{mylemma}
\label{lemma:Qpkodd}
For any $p$ prime and $k$ odd we have
\begin{equation*}
Q(p^k) =0.
\end{equation*}
\end{mylemma}
We follow the argument of \cite{Y2}, Section 9.2.
\begin{proof}
We may assume $p$ is odd.  Using the Hecke relations, we see that
\begin{equation*}
\lambda_{a, b}(p^{2k + 1}) = \sum_{j \leq k} c_j (\lambda_{a, b}(p))^{2j + 1}
\end{equation*}
for some complex numbers $c_j$ not depending on $a$ and $b$, provided $(4a^3 + 27b^2, p) = 1$ (it is a Chebyshev polynomial of the second kind).  
Hence
\begin{equation*}
Q(p^{2k + 1}) = \sum_{j \leq k} c_j \mathop{\sum_{\alpha \shortmod{p}}\sum_{\beta \shortmod{p}}}_{(4\alpha^3 + 27 \beta^2, p) = 1} (\lambda_{\alpha, \beta}(p))^{2j + 1} + \mathop{\sum_{\alpha \shortmod{p}}\sum_{\beta \shortmod{p}}}_{4\alpha^3 + 27 \beta^2 \equiv 0 \shortmod{p}} (\lambda_{\alpha, \beta}(p))^{2k + 1}.
\end{equation*}
%It therefore suffices to show that
%\begin{equation*}
%\sum_{\alpha \shortmod{p}}\sum_{\beta \shortmod{p}} (\lambda_{\alpha, \beta}(p))^{2k + 1} = 0.
%\end{equation*}
%Opening the character sum formula for $\lambda(p)$, we see that the above sum is
%\begin{equation*}
%\left(\frac{-1}{\sqrt{p}} \right)^{2k + 1} \mathop{\sum \cdots \sum}_{x_1, \ldots, x_{2k + 1}} \sum_{\alpha %\shortmod{p}}\sum_{\beta \shortmod{p}} \prod_{j=1}^{2k + 1} \left(\frac{x_i^3 + \alpha x_i + \beta}{p} \right).
%\end{equation*}
%Let $e$ be a quadratic nonresidue $\pmod{p}$ and apply the changes of variables $x_i \rightarrow e x_i$, $\alpha \rightarrow e^2 \alpha$, $\beta \rightarrow e^3 \beta$ to the above sum.  We obtain the negative of the same sum so it vanishes. 
The latter sum vanishes by the parameterization given by Lemma \ref{lemma:parameterization} and by using a simple generalization of \eqref{eq:degeneratecalculation}.  The first sum vanishes because elliptic curves modulo $p$ come in pairs with opposite values of $\lambda(p)$ (which can be seen by quadratic twisting).
\end{proof}

\subsection{Summing over $a$ and $b$}
At this point we develop the summation over $a$ and $b$ present in \eqref{eq:LUMapprox}, namely
\begin{equation*}
\mathop{\sum \sum}_{\substack{(a,b) \in S \\ \Delta \equiv 0 \shortmod{c}}} \lambda_{a,b}(mn) w_X(a,b) = \sum_{(g,2mn) = 1} \mu(g) \mathop{\sum_a \sum_{b \text{ odd }} }_{g^6 \Delta \equiv 0 \shortmod{c}} \lambda_{ag^2, bg^3}(mn) w_X(ag^2,bg^3).
\end{equation*}
Note that we can make the restriction $g \leq X^{\varepsilon}$ without changing the error term in Lemma \ref{lem:LUMavgclean}.

\begin{mydefi}
Let $\mathcal{C}_x$ be the set of all Dirichlet characters of moduli $\leq x$ and let $\mathcal{N}_x = \{n \in \mz : |n| \leq x \}$.  Let $\mathcal{Y}^{(i, j)}_x = \mathcal{N}_x^{i} \times \mathcal{C}_x^{j}$ (here the superscripts on $\mathcal{N}_x$ and $\mathcal{C}_x$ indicate cross product).  
\end{mydefi}

We use the following
\begin{myprop}
\label{prop:MTRT}
Let $r$ and $c$ be positive odd integers such that $c \leq X^{\varepsilon}$ is squarefree.  Suppose $(g, 2r) = 1$ and $g \leq X^{\varepsilon}$. Set
\begin{equation*}
Z = \mathop{\sum_a \sum_{b \text{ odd }}}_{g^6 \Delta \equiv 0 \shortmod{c}} \lambda_{ag^2, bg^3}(r) w_X(ag^2, bg^3) %\left(\frac{ag^2}{A}, \frac{bg^3}{B} \right).
\end{equation*}
Let $r = r_1 r_2$ where $r_1$ is the product of primes exactly dividing $r$.  Set $c = c_0 c_1 c_2 (c,g)$, where $c_1 | r_1$, $c_2 | r_2$, and $(r, c_0) = 1$.  Further split $r_1 = c_1 r_0$.  There exists a set $\mathcal{Y}$ of the form $\mathcal{Y}^{(i, j)}_{x}$ ($i = 9^4, j= 3^4$ is sufficient) with $x \ll_{\varepsilon} r_2^{2004} X^{4\varepsilon}$ and a function $F$ on $\mathcal{Y}$ satisfying
\begin{equation*}
|F(y)| \ll r_2^{2004} X^{\varepsilon}
\end{equation*}
for $y \in \mathcal{Y}$,
%\begin{equation*}
%|\mathcal{Y}| \ll |r_2|^{2004} X^{\varepsilon}.
%\end{equation*}
%Further, the integers in $\mathcal{Y}$ and the moduli of the characters in $\mathcal{Y}$ are $\ll |r_2|^{2004} X^{\varepsilon}$.  
such that
\begin{equation*}
Z= M.T. + R.T + O(r_2X^{1/2}),
\end{equation*}
where
\begin{equation}
\label{eq:MT}
M.T. = \onehalf AB \widehat{w}(0,0) \delta(r_1 = 1) \frac{Q_{c}(r)}{c_0 r^*{}^2}  \frac{1}{g^5}
\end{equation}
and
\begin{multline*}
R.T. = %\frac{AB}{r_0} \varepsilon_{r_0} \mu(r_0) \sum_{y \in \mathcal{Y}} F(y) \sum_{\substack{k_1 \in \mz \\ k_1 \neq 0}} \sum_{e | k_1^2 m_{1, y}} \chi_{3, y}(e) \chi_{4, y}(t(e)) \sum_{\substack{h_1 \in \mz \\ h_1 \neq 0}} \\
 %\sum_{\chi \shortmod{(k_1^2/e) m_{2, y}}} \tau(\chi) \left(\frac{k_1}{r_0} \right) \chi(r_0) \chi_{1, y}(r_0) \overline{\chi}^3(h_1) \chi_{2, y}(h_1) w_2(h_1, k_1, r_0, y),
\frac{AB}{r_0} \varepsilon_{r_0} \mu(r_0) 
\sum_{y \in \mathcal{Y}} F(y) \sum_{k_1} \chi_{1,y}(k_1) \left(\frac{k_1}{r_0}\right) \sum_{e | k_1^2 m_{1,y} } \chi_{2,y}(t(e)) \chi_{3,y}(e) \sum_{h_1} \chi_{4,y}(h_1) \chi_{5,y} (r_0)  
\\
\cdot \frac{1}{\varphi(m_{2,y} k_1^2  /e)}\sum_{\chi \shortmod{m_{2,y} k_1^2/e}} \tau(\chi) \overline{\chi}^3(h_1) \chi(r_0) \overline{\chi}(t^3(e)/e).
\end{multline*}
and where
\begin{equation*}
w_2 = e\left(\frac{h_1^3 t^3(e) m_{3,y}}{k_1^2 r_0 m_{4,y}}\right) \widehat{w}\left(\frac{h_1A m_{5,y}}{r_0 m_{6,y}}, \frac{k_1B m_{7,y}}{r_0 m_{8,y}} \right).
\end{equation*}
Here $\chi_{i, y}, m_{i, y}$ are components of vectors in $\mathcal{Y}$ and $t(e)$ is the least positive integer $l$ such that $e | l^3$.  The implied constant in \eqref{eq:MT} depends on $w$ only.
The function $F$ is allowed to depend on $c_1$, $c_2$, $c_3$, $g$, and $r_2$ but not $r_0$.
\end{myprop}
Remarks.  In our applications of Proposition \ref{prop:MTRT} we shall have $|r_2| \ll X^{\varepsilon}$ so that the expansion into characters has little cost.  

The form of the remainder term above is well-suited for averaging over $r$ because the variables have been separated.

The set $\mathcal{Y}$ and function $F$ could be made explicit but the size of the formulas would be prohibitive. 
%In our application of this result we shall obtain savings in the summations over $r_1$, $h_1$, and $k_1$ only; the bounds will be uniform with respect to $y \in \mathcal{Y}$.

%Without changing the form of $R.T.$ we may replace $r_0$ with $r_1$; this simply changes $F$, $m_{i,y}$, and $\chi_{i,y}$ slightly.

We record here two formulas that will be useful in the proof.  The following elementary reciprocity law is valid for coprime integers $u$ and $v$
\begin{equation}
\label{eq:elementaryreciprocity}
\frac{\overline{u}}{v} + \frac{\overline{v}}{u} \equiv \frac{1}{uv} \pmod{1}.
\end{equation}

We also make extensive use of
\begin{equation}
\label{eq:addtomult}
e \left(\frac{a}{n} \right) = \frac{1}{\varphi(n)} \sum_{\chi \shortmod{n}} \tau(\chi) \overline{\chi}(a),
\end{equation}
valid for $(a, n) = 1$.  Here $\tau(\chi)$ is the Gauss sum.  This formula is useful for separating the variables present in an exponential, as long as $n$ is not too large.  In our applications it becomes a notational nuisance to ensure $(a, n) = 1$ (by dividing out by greatest common divisors), but the benefit of separating the variables is well worth the price.

\begin{proof}
We begin by setting $c = (c, g) c_4$ with $c_4 = c_0 c_1 c_2$.  We obtain
\begin{equation*}
Z = \mathop{\sum_a \sum_{b \text{ odd}}}_{ \Delta \equiv 0 \shortmod{c_4}} \lambda_{ag^2, bg^3}(r) w_X(ag^2,bg^3).
\end{equation*}
Applying Poisson summation in $a$ modulo $q= r_1 r_2^* c_0$ and $b$ modulo $2q$ leads to the identity (recall $n^*$ denotes the product of primes dividing $n$)
\begin{equation}
\label{eq:Zpoisson}
Z = \frac{AB}{2g^5 q^2} \sum_h \sum_k \mathop{\sum_{\alpha \shortmod{q}} \sum_{\beta \shortmod{2q}}}_{\substack{ \Delta \equiv 0 \shortmod{c_4} \\ \beta \text{ odd}}} \lambda_{\alpha g^2, \beta g^3}(r)
e \left(\frac{2\alpha h + \beta k}{2q} \right)
 \widehat{w}\left(\frac{hA}{g^2q}, \frac{kB}{2g^3q} \right);
\end{equation}
here $\widehat{w}$ is the Fourier transform of $w$.
We separate the variables by letting $\alpha = \alpha_1 (q/r_1)(\overline{q/r_1}) + \alpha_2 (q/r_2^*)(\overline{q/r_2^*}) + \alpha_3 (q/c_0)(\overline{q/c_0})$, and similarly for $\beta$ but with an extra summand corresponding to $2$.  Then
\begin{equation*}
Z = \frac{AB}{2g^5 q^2} \sum_h \sum_k S_1 S_2 S_3 e\left(\frac{k}{2} \right) \widehat{w}\left(\frac{hA}{g^2q}, \frac{kB}{2g^3q} \right),
\end{equation*}
where
\begin{equation*}
S_1 = \mathop{\sum_{\alpha_1 \shortmod{r_1}} \sum_{\beta_1 \shortmod{r_1}}}_{\Delta \equiv 0 \shortmod{c_1}}
\lambda_{\alpha_1 g^2, \beta_1 g^3}(r_1) e \left(\frac{(\alpha_1 h + \bar{2} \beta_1 k)\overline{c_0 r_2^*}}{r_1} \right),
\end{equation*}
\begin{equation*}
S_2 = \mathop{\sum_{\alpha_2 \shortmod{r_2^*}} \sum_{\beta_2 \shortmod{r_2^*}}}_{\Delta \equiv 0 \shortmod{c_2}}
\lambda_{\alpha_2 g^2, \beta_2 g^3}(r_2) e \left(\frac{(\alpha_2 h + \bar{2} \beta_2 k)\overline{r_1 c_0}}{r_2^*} \right),
\end{equation*}
and
\begin{equation*}
S_3 = \mathop{\sum_{\alpha_3 \shortmod{c_0}} \sum_{\beta_3 \shortmod{c_0}}}_{\Delta \equiv 0 \shortmod{c_0}}
 e \left(\frac{(\alpha_3 h + \bar{2} \beta_3 k)\overline{r_1 r_2^*}}{c_0} \right).
\end{equation*}
We use Corollary \ref{cor:maincompletesum} to compute $S_1$, Lemma \ref{lem:degenerateexpsum} to compute $S_3$, and leave $S_2$ as it is (since we do not have a closed formula for $S_2$).  We obtain
\begin{multline}
\label{eq:usefulZ}
Z = \frac{AB \varepsilon_{r_0} \mu(r_0) r_0}{2g^5 q^2 \sqrt{c_1}} \frac{}{} \sum_h \sum_k
\left(\frac{3g}{c_1} \right) \left(\frac{2gkc_0 c_1 r_2^*}{r_0} \right)
e\left( \frac{-4h^3 \overline{k}^2 \overline{c_0c_1 r_2^* }}{r_0} \right) S_2  e\left(\frac{k}{2}\right)
\widehat{w}\left(\frac{hA}{g^2q}, \frac{kB}{2g^3q} \right)\\
\cdot
\sum_{\gamma_1 \shortmod{c_1}} \left(\frac{\gamma_1}{c_1} \right)
e\left( \frac{(-3\gamma_1^2h + \gamma_1^3 k)\overline{c_0 r_0 r_2^* }}{c_1} \right) \sum_{\gamma_0 \shortmod{c_0}} e \left(\frac{(-3\gamma_0^2 h  + \gamma_0^3 k )\overline{r_1 r_2^*} }{c_0} \right).
\end{multline}
The zero frequency (the term with $h = k = 0$) gives the main term
\begin{equation}
\label{eq:mainterm}
\onehalf AB \widehat{w}(0,0) \delta(r_1 = 1) \frac{Q_{c}(r)}{c_0r^*{}^2} \frac{1}{g^5},
\end{equation}
using that $\left(\frac{0}{r_0} \right) = \delta(r_0 = 1)$, that the summation of $\gamma_1$ is zero unless $c_1 = 1$ (hence $r_1 =1$), and that $S_2 = Q_{c_2}(r_2)$ (which equals $Q_c(r)$ since $r=r_2$ and $c_2 = (c, r_2)$).

By trivial estimates (bounding all sums by their absolute values) the terms with $h = 0$ and $k \neq 0$ contribute at most $O(r_2^*A )$, with an implied constant depending on $w$ only.  Similarly, the terms with $k= 0$ and $h \neq 0$ contribute $O(r_2^*B )$.  An improvement on the trivial bound $|S_2| \leq r_2^*{} ^2$ would give savings in these two bounds but in our applications $r_2^* \ll X^{\varepsilon}$ so these bounds are more than sufficient.

Let $Z_1$ be the terms of $Z$ with $h \neq 0$ and $k \neq 0$; we continue with $Z_1$.  Our goal is twofold.  Firstly, we would like to use \eqref{eq:addtomult} in order to write the various exponential sums modulo $c_i$ in terms of multiplicative characters.  The difficulty here is that the numerator and denominator in the exponentials may not be coprime.  We must introduce a lot of notation to place the fractions into lowest terms.  The second device is to use \eqref{eq:elementaryreciprocity} to write
\begin{equation}
\label{eq:flip}
e\left( \frac{-4h^3 \overline{k}^2 \overline{c_0c_1 r_2^* }}{r_0} \right) = e\left(\frac{4h^3 \overline{r_0}}{k^2 c_0 c_1 r_2^*}\right) e\left(\frac{-4h^3 }{k^2 c_0 c_1 r_0 r_2^*}\right),
\end{equation}
a device that effectively reduces the size of the modulus.  Only then do we expand this exponential factor into Dirichlet characters, but again there is a coprimality issue.

As an initial simplification, we sum over $k$ even and $k$ odd separately; set $Z_1 = Z_2 + Z_3$, where $Z_2$ corresponds to $k$ even and $Z_3$ corresponds to $k$ odd.  As a slight notational simplification we shall treat $Z_2$ only; the sum $Z_3$ can be treated similarly. We first change variables by $k \rightarrow 2k$.  Set $e = (h^3, k^2)$, and for any integer $n$ let $t(n)$ be the least positive integer $n_0$ such that $n | n_0^3$ (so for any $m$, $n$, if $n | m^3$ then $t(n) | m$).
Since $t(e) | h$ we may set $h = t(e) h_0$.  The condition $e = (h^3, k^2)$ is seen to be equivalent to the two conditions $(h_0, k^2/e) = 1$ and $(t^3(e)/e, k^2/e) = 1$.  Applying \eqref{eq:flip}, we have
\begin{multline*}
Z_2 = \frac{AB}{2g^5 q^2} \varepsilon_{r_0} \mu(r_0) \frac{r_0}{\sqrt{c_1}} \left(\frac{3g}{c_1} \right) \mathop{\sum_k \sum_{e | k^2}}_{(t^3(e)/e, \, k^2/e) = 1} \sum_{(h_0, k^2/e) = 1}
 \left(\frac{gkc_0 r_2^*}{r_0} \right)
e\left(\frac{h_0^3 (t^3(e)/e) \overline{r_0}}{(k^2/e) c_0 c_1 r_2^*}\right) S_2\\
\cdot
\sum_{\gamma_1 \shortmod{c_1}} \left(\frac{\gamma_1}{c_1} \right)
e\left( \frac{(-3\gamma_1^2h_0 t(e) + 2\gamma_1^3 k)\overline{c_0 r_0 r_2^* }}{c_1} \right) 
\sum_{\gamma_0 \shortmod{c_0}} e \left(\frac{(-3\gamma_0^2 h_0 t(e)  + 2\gamma_0^3 k )\overline{r_1 r_2^*} }{c_0} \right)  w_1,
\end{multline*}
where
\begin{equation*}
w_1 = e\left(\frac{h_0^3 t^3(e) }{k^2 c_0 c_1 r_0 r_2^*}\right) \widehat{w}\left(\frac{h_0 t(e)A}{g^2q}, \frac{kB}{g^3q} \right).
\end{equation*}
Now set $e_0 = (h_0^3, c_0)$, $e_1 = (h_0^3, c_1)$, $e_2 = (h_0^3, r_2^*)$ (note $e_i | h_0$ since each $e_i$ is squarefree) and write $h_0 = e_0 e_1 e_2 h_1$ (which can be done because the $e_i$ are pairwise coprime).  These changes of variables transform $Z_2$ to (we omit the display of the summands)
\begin{equation*}
\sum_{\substack{e_0 | c_0 \\ e_1 | c_1 \\ e_2 | r_2^*}} \mathop{\sum_k \sum_{e | k^2}}_{(\frac{k^2}{e}, \frac{t^3(e)}{e} e_0 e_1 e_2) = 1} \sum_{(h_1, \frac{k^2}{e} \frac{c_0}{e_0}\frac{c_1}{e_1}\frac{r_2^*}{e_2} ) = 1}.
\end{equation*}
Similarly, set $f_0 = (k, c_0)$, $f_1 = (k, c_1)$, $f_2 = (k, r_2^*)$, and
write $k = f_0 f_1 f_2 k_1$.  Then the summation conditions become
\begin{equation}
\label{eq:coprimalityconditions}
\sum_{\substack{e_0 | c_0 \\ e_1 | c_1 \\ e_2 | r_2^*}} \sum_{\substack{f_0 | c_0 \\ f_1 | c_1 \\ f_2 | r_2^*}} \mathop{\sum_{(k_1, \frac{c_0}{f_0}\frac{c_1}{f_1}\frac{r_2^*}{f_2}) = 1} \sum_{e | k_1^2 f_0^2 f_1^2 f_2^2}}_{(\frac{k_1^2 f_0^2 f_1^2 f_2^2}{e}, \frac{t^3(e)}{e} e_0 e_1 e_2) = 1} \sum_{(h_1, \frac{k_1^2 f_0^2 f_1^2 f_2^2}{e} \frac{c_0}{e_0}\frac{c_1}{e_1}\frac{r_2^*}{e_2} ) = 1}.
\end{equation}
The point is that these new variables are sufficiently coprime that we can expand the exponentials into multiplicative characters.

We now focus on the various summands of $Z_2$.  We have
\begin{multline*}
\sum_{\gamma_1 \shortmod{c_1}} \left(\frac{\gamma_1}{c_1} \right)
e\left( \frac{(-3\gamma_1^2h_0 t(e) + 2\gamma_1^3 k)\overline{c_0 r_0 r_2^* }}{c_1} \right) \\
=
\sum_{\gamma_1 \shortmod{c_1}} \left(\frac{\gamma_1}{c_1} \right) e\left( \frac{(-3\gamma_1^2 h_1 t(e) e_0 e_2 )\overline{c_0 r_0 r_2^* }}{c_1/e_1} \right) e\left( \frac{( 2\gamma_1^3 f_0 f_2 k_1)\overline{c_0 r_0 r_2^* }}{c_1/f_1} \right).
\end{multline*}
It is readily apparent from \eqref{eq:usefulZ} that the sum vanishes unless $(e, c_1) = 1$.
Because of the coprimality conditions of \eqref{eq:coprimalityconditions} we may apply \eqref{eq:addtomult} to each of the above exponentials, and obtain
\begin{multline*}
\frac{1}{\varphi \left(\frac{c_1}{e_1}\right) \varphi \left( \frac{c_1}{f_1}\right)} \sum_{\chi_1 \shortmod{\frac{c_1}{e_1}}} \sum_{\chi_2 \shortmod{\frac{c_1}{f_1}}}
\tau(\chi_1) \tau(\chi_2) \\
\cdot
\sum_{\gamma_1 \shortmod{c_1}} \left(\frac{\gamma_1}{c_1} \right) \overline{\chi_1}(-3\gamma_1^2 h_1 t(e) e_0 e_2\overline{c_0 r_0 r_2^* }) \overline{\chi_2}((2\gamma_1^3 f_0 f_2 k_1)\overline{c_0 r_0 r_2^*}).
\end{multline*}
This sum has the form
\begin{equation}
\label{eq:sumform}
\sum_{y_1 \in \mathcal{Y}_1} F_1(y_1) (\chi_1 \chi_2) (r_0) \overline{\chi_1}(h_1 t(e)) \overline{\chi_2}(k_1),
\end{equation}
where $\mathcal{Y}_1$ and $F_1$ satisfy the same conditions as $\mathcal{Y}$ and $F$ stated in the Proposition (with $i=9, j=3$, and $x= r_2 X^{2\varepsilon}$).  Precisely, for 
\begin{equation*}
y_1 = (e_0, e_1, e_2, \gamma_1, c_0, c_1, r_2^*, f_0, f_2, \chi_1, \chi_2, \chi_1 \chi_2)
\end{equation*}
 we have
 \begin{equation*}
F_1(y_1) = 
%\frac{1}{\varphi(c_1/e_1) \varphi(c_1/f_1)}  
\frac{1}{\varphi \left(\frac{c_1}{e_1}\right) \varphi \left( \frac{c_1}{f_1}\right)}
\tau(\chi_1) \tau(\chi_2) \left(\frac{\gamma_1}{c_1}\right) \overline{\chi_1}(-3\gamma_1^2 e_0 e_2 \overline{c_0 r_2^* }) \overline{\chi_2}(2\gamma_1^3 f_0 f_2)\overline{c_0  r_2^*}).
\end{equation*}
The set $\mathcal{Y}_1$ can be taken to be all such vectors $y_1$ above if we adjust $F_1$ to be zero unless the restrictions given by \eqref{eq:coprimalityconditions} hold, as well as other obvious restrictions such as that that $\chi_1$ is modulo $c_1/e_1$, $\gamma_1$ runs modulo $c_1$, etc.  In addition, we should let the final component of $y_1$ be an arbitrary character $\chi_3$ and have $F$ detect $\chi_3 = \chi_1 \chi_2$.

Using similar reasoning, we may write
\begin{equation*}
\sum_{\gamma_0 \shortmod{c_0}} e \left(\frac{(-3\gamma_0^2 h_0 t(e)  + 2\gamma_0^3 k )\overline{r_1 r_2^*} }{c_0} \right)
= \sum_{y \in \mathcal{Y}_2} F_2(y) (\chi_3 \chi_4) (r_0) \overline{\chi_3}(h_1 t(e)) \overline{\chi_4}(k_1).
\end{equation*}

The sum $S_2$ is slightly more complicated to write explicitly because there is the additional notational issue of the greatest common factors of $\alpha_2$ and $r_2^*$ and likewise with $\beta_2$.  Nevertheless, it is not difficult to see that it can be written in the same form as \eqref{eq:sumform}.

Finally, we have
\begin{align*}
e\left(\frac{h_0^3 (t^3(e)/e) \overline{r_0}}{(k^2/e) c_0 c_1 r_2^*}\right) = e\left(\frac{h_1^3 (t^3(e)/e)(e_0^2 e_1^2 e_2^2)\overline{r_0}}{(k^2/e) (c_0/e_0) (c_1/e_1) (r_2^*/e_2)}\right) \\
= \frac{1}{\varphi(u)} \sum_{\chi \shortmod{u}} \tau(\chi) \overline{\chi}(h_1^3 \overline{r_0} (t^3(e)/e)) F_3(c_0, c_1, c_2, e_0, e_1, e_2, f_0, f_1, f_2),
\end{align*}
where $ u = (k_1^2 f_1^2 f_2^2 f_3^2/e)  (c_0/e_0) (c_1/e_1) (c_2/e_2)$ and $F_3$ is a function satisying the conditions of the Proposition.  Here we cannot absorb $u$, $\chi \pmod{u}$, and $\tau(\chi)$ into $\mathcal{Y}$ because $k_1$ can grow larger than $X^{\varepsilon}$.

Now we combine the above formulas with the summation conditions \eqref{eq:coprimalityconditions}.  We obtain a sum of the type
\begin{multline*}
\sum_{y \in \mathcal{Y}} F(y) \sum_{k_1} \chi_{1,y}(k_1) \sum_{e | k_1^2 m_{1,y} } \chi_{2,y}(t(e)) \chi_{3,y}(e) \sum_{h_1} \chi_{4,y}(h_1) \chi_{5,y} (r_0)  
\\
\frac{1}{\varphi(m_{2,y} k_1^2  /e)}\sum_{\chi \shortmod{m_{2,y} k_1^2/e}} \tau(\chi) \overline{\chi}^3(h_1) \chi(r_0) \overline{\chi}(t^3(e)/e).
\end{multline*}
The various characters $\chi_{i, y}$ arise as products of a bounded number of characters of moduli $\ll r_2^* X^{\varepsilon}$.  The integers $m_{i,y}$ are products of a bounded number of $e_i$'s and $f_i$'s, etc.
We may assume further that the coprimality restrictions listed in \eqref{eq:coprimalityconditions} hold either by use of $F$ (to handle coprimality amongst $e_i$'s, $f_i$'s, etc.) or by the use of the characters $\chi_{i,y}$ (to handle coprimality of $k_1$ and $h_1$ with respect to integers coming from $y \in \mathcal{Y}$).  Further, we write $q^2 = r_0^2 c_0^2 c_1^2 r_2^*{}^2$ and absorb the factor $c_0^2 c_1^2 r_2^*{}^2$ into $F$; likewise, we absorb $(2 g^5 \sqrt{c_1})^{-1}$ into $F$.  The quadratic character $\left(\frac{gkc_0r_2^*}{r_0}\right)$ is of the form $\left(\frac{k_1}{r_0}\right) \chi_{y}(r_0)$. 

Similarly, an inspection of how $w_1$ is altered under the changes of variable replacing $h_0$ with $h_1$ and $k$ with $k_1$ leads to $w_2$ being of the stated form; it is important to notice the exponential factor
\begin{equation*}
e \left(\frac{h_1^3 t^3(e) m_{3,y}}{k_1^2 r_0 m_{4,y}} \right)
\end{equation*}
is absorbed into $w_2$ also.  The proof is now complete.
%In the remarks after the statement of the Proposition we mentioned that we may replace $r_0$ with $r_1$ without changing the statement.  This follows because we allowed $F$ and the characters $\chi_{i, y}$ to depend on $c_1$; by using standard techniques of separation of variables the transition from $r_0$ to $r_1$ can easily be accomplished.
\end{proof}

\subsection{Estimating the remainder term}
\label{section:RT}
We combine the result of Proposition \ref{prop:MTRT} with Lemma \ref{lem:LUMavgclean}.  The main term is evaluated in the following section.  The $O(r_2X^{1/2})$ term is clearly small.  We estimate the contribution from $R.T.$, say $R$, here.  We apply Proposition \ref{prop:MTRT} with $r_1 = (mn)_1$ and $r_2 = (mn)_2$.  We let $(mn)_1 = ((mn)_1, c) (mn)_0$. %(using the remark stated after Proposition \ref{prop:MTRT} saying that we may freely replace $r_0$ with $r_1$ without altering the form of $R.T.$  
We have
\begin{multline}
\label{eq:remainder}
R = AB \sum_{d} \sum_{l} \sum_{\substack{c | dl \\ c \leq X^{\varepsilon}}} \mu(c) \sum_{\substack{g \leq X^{\varepsilon} \\ (g, 2) = 1}} \mu(g) \mathop{\sum_{m} \sum_{n}}_{(mn, g) = 1} \frac{\mu^2(dlm) \mu(dm)}{dl\sqrt{mn}(mn)_0}  \varepsilon_{(mn)_0} \mu((mn)_0) 
 \sum_{y \in \mathcal{Y}} f(y) \sum_{k_1 \neq 0}
\\
 \sum_{e | k_1^2 m_{1,y}} \chi_{3, y}(e)  
 \chi_{4, y}(t(e)) \sum_{h_1} 
 \sum_{\chi \; (m_{2, y} k_1^2/e )} \tau(\chi) \left(\frac{k_1}{(mn)_0} \right) (\chi \chi_{1, y})((mn)_0) (\overline{\chi}^3\chi_{2, y})(h_1) \chi(t^3(e)/e) w_3,
\end{multline}
where
\begin{equation*}
w_3 = e\left(\frac{h_1^3 m_{3,y}}{k_1^2 (mn)_0 m_{4,y} }\right) \widehat{w}\left(\frac{h_1A m_{5,y}}{(mn)_0 m_{6,y}}, \frac{k_1B m_{7,y}}{(mn)_0 m_{8,y}} \right) Y\left(\frac{2\pi dn}{U} \right) P \left( \frac{\log{M/dml^2}}{\log{M}} \right).
\end{equation*}
The basic idea at this point is to use the Riemann hypothesis for Dirichlet L-functions to show that there is square-root cancellation in the summation over $h_1$ and $(mn)_0$ and treat the summation over the other variables trivially (of course the only other variables that could contribute anything substantial are $k_1$ and $\chi$).  The details are rather technical and have essentially already been carried out in \cite{Y}.  There we estimated a sum very similar to \eqref{eq:remainder}; the sum (19) in \cite{Y} had a sum over $p$ prime instead of $(mn)_0$ but is otherwise completely analogous to \eqref{eq:remainder} here.  On a superficial level there are two essential differences between the two sums.  The transition from a sum over primes to squarefree integers poses no difficulty; since we used L-function methods to bound the sum over primes in \cite{Y} it is simple to alter the proofs to suit squarefree $n$ (instead of dealing with contour integrals of $L'/L(s, \chi)$ we have $1/L(s, \chi)$).  Any potential difficulty arising from the fact that the sum in \eqref{eq:remainder} is over the product $(mn)_0$ is easily handled by taking products of L-functions.  We sketch the arguments used in \cite{Y}; for full details see Section 5.2 of the original paper.

We separate our arguments into four cases, depending on whether one or both of the characters $\psi(n) = (k_1/n) \chi(n) \chi_{1, y}(n)$ and $\overline{\chi}^3 \chi_{2, y}$ are principal.  For the typical case where neither of the two characters is principal we use the Riemann hypothesis to get square-root cancellation in the summation over $m$, $n$, and $h_1$.  One difficulty is the presence of the exponential in the test function $w_3$.  The exponential factor can be effectively controlled using stationary phase estimates.  It is here that the restriction $UM \ll X^{7/9 - \varepsilon}$ is the limit of the method.

The case where $\overline{\chi}^3 \chi_{2, y}$ is principal but $\psi$ is not is easily handled by trivial estimations.

Consider the case where both characters are principal.  Then $\chi$ effectively has modulus $\ll X^{\varepsilon}$ and $k_1 = \square m_y$, where $m_y \ll X^{\varepsilon}$.  The bound obtained by this method is sufficient only in certain ranges, namely for $k_1 \gg UMX^{-2/3 + \varepsilon}$.  To handle $k_1$ small we use a simple application of Weyl's method to estimate sums of the type
\begin{equation*}
\sum_{N \leq n \leq 2N} \left| \sum_{H \leq h \leq 2H} e\left(\frac{h^3 \overline{k}^2}{n} \right)\right|.
\end{equation*}

The remaining case where $\psi$ is principal but the other character is not principal is also handled exactly as in \cite{Y}.  It is necessary to obtain extra saving arising from the oscillation of the exponential factor
\begin{equation*}
e\left(\frac{h_1^3 m_{1,y}}{k_1^2 (mn)_1 }\right).
\end{equation*}
We can directly apply Lemma 5.5 of \cite{Y} to bound the necessary sum in this situation.

To conclude, we have sketched a proof that $|R| \ll AB X^{-\delta}$.  The full details of the proof are exceedingly similar to those in \cite{Y} and can be repeated {\em mutatis mutandi}.

\subsection{Evaluating the main term}
\label{section:mainterm}
In this section we compute the main term, namely the quantity formed by inserting the $M.T.$ into \eqref{eq:LUMapprox}.
This quantity is
\begin{equation}
\label{eq:mollifiedmainterm}
\frac{AB}{2} \widehat{w}(0,0) \sum_{d} \sum_{l} \sum_{\substack{c | dl \\ c \leq X^{\varepsilon}}} \mu(c) \sum_{(g, 2) =1} \frac{\mu(g)}{g^5} \mathop{\sum_{m} \sum_{n}}_{\substack{(mn)_1 = 1 %\\ (mn)_2 \ll X^{\varepsilon} \\ 
(m, dl) = 1 \\ (mn, g) = 1}} \frac{ \mu^2(dlm) \mu(dm)}{dl\sqrt{mn}} \frac{Q_c(mn)}{c_0 (mn)^*{}^2} H_1(d, l, m, n),
\end{equation}
where $H_1(d,l,m,n)$ is given as in Lemma \ref{lem:LUMapprox}.
%\begin{equation*}
%H_1(d, l, m, n) = Y\left(\frac{2\pi dn}{U} \right) P \left( \frac{\log{M/dml^2}}{\log{M}} \right).
%\end{equation*}
As a first step towards simplification we execute the summation over $c$.  We compute
\begin{equation*}
\sum_{c | dl} \frac{\mu(c)}{c_0} Q_c(mn) = \sum_{\alpha} \sum_{\beta} \lambda_{\alpha,\beta}(mn) \sum_{\substack{c | dl \\ (c, mn) | D}} \frac{\mu(c)}{c_0},
\end{equation*}
where $D = 4\alpha^3 + 27\beta^2$.  %Now let $dl = (dl, mn) u$ and write $c = c_{mn} c_u $ where $c_{mn} | (dl, mn)$ and $(c_u, mn) = 1$.  Furthermore, notice that $c_0$ does not depend on $c_{mn}$ because $(c_0, mn) = 1$.
Set $c_g = (c, g)$ (so that $c = c_0 c_2 c_g$ since $c_1 = 1$) and $dl = (dl, mn) (dl, g) \frac{dl}{(dl, mng)}$, so that the above sum factors as
\begin{equation*}
 \sum_{\alpha} \sum_{\beta} \lambda_{\alpha,\beta}(mn) \sum_{\substack{c_0 | dl \\ (c_0, mng) = 1}} \frac{\mu(c_0)}{c_0} \sum_{\substack{c_2 | (dl, mn) \\ c_2 | D}} \mu(c_2) \sum_{c_g | (dl, g)} \mu(c_g).
\end{equation*}
Define
\begin{equation*}
Q_k'(r) = \sum_{\alpha \shortmod{r^*}} \sum_{\beta \shortmod{r^*}} \psi_{\Delta}((k, r)) \lambda_{\alpha, \beta}(r),
\end{equation*}
and let
\begin{equation*}
\varphi^*(n) = \frac{\varphi(n)}{n}.
\end{equation*}
Then the desired sum over $c$ is
\begin{equation*}
Q_{dl}'(mn) \varphi^* \left(\frac{dl}{(dl, mn)} \right) \delta((dl, g) = 1).
\end{equation*}
Note that $Q_k'(r)$ is multiplicative in $r$ for any fixed $k$ and that $Q_1'(r) = Q(r)$.

Inserting this expression into \eqref{eq:mollifiedmainterm} (after extending $c > X^{\varepsilon}$ with $O(ABX^{-\delta})$ error), we obtain that the main term is
\begin{multline}
\label{eq:mollifiedmaintermtwo}
\frac{AB}{2}  \widehat{w}(0,0) \sum_{d} \sum_{l} \sum_{(g, 2dl) =1} \frac{\mu(g)}{g^5} \mathop{\sum_{m} \sum_{n}}_{\substack{(mn)_1 = 1 \\ %(mn)_2 \ll X^{\varepsilon} \\ 
(gm, dl) = 1 \\ (mn, g) = 1}} \frac{ \mu^2(dlm) \mu(dm) \varphi^* \left(\frac{dl}{(dl, mn)} \right)}{dl\sqrt{mn}}  \frac{Q'_{dl}(mn)}{(mn)^*{}^2} H_1(d, l, m, n) \\
+ O(ABX^{-\delta}).
\end{multline}

We may freely remove the truncation factor $\phi((mn)_2/X^{\varepsilon})$ from $H_1(d,l,m,n)$ without introducing a new error term; that is, we may replace $H_1$ with $H$, where
\begin{equation*}
H(d,l,m,n) = Y\left(\frac{2\pi dn}{U} \right) P \left( \frac{\log{M/dml^2}}{\log{M}} \right)
\end{equation*}
without introducing a new error term in \eqref{eq:mollifiedmaintermtwo}.  Let $\Theta$ be this sum.  Recall $P(x)$ for $x \geq 0$ is a polynomial satisfying $P(0) = 0$, $P(1) = 1$ and that $P(x)$ is extended to be zero for $x < 0$.  Using the Mellin transform of $H$ we have
\begin{equation*}
\Theta = \frac{AB}{2}  \widehat{w}(0,0) \sum_{j \geq 1} \frac{a_j j!}{(\log{M})^j} \frac{1}{(2 \pi i)^2} \int_{(1)} \int_{(1)} Z(s,v) \frac{M^s}{s^{j +1}} \frac{U^v}{v} (2\pi)^{-v} \Gamma(1 + v) G(v) ds dv,
\end{equation*} 
where
\begin{equation*}
Z(s,v) = \sum_{d} \sum_{l} \sum_{(g, 2dl) =1} \mathop{\sum_{m} \sum_{n}}_{\substack{(mn)_1 = 1  \\  (mn, g) = 1}} \frac{ \mu^2(dlm) \mu(dm) \varphi^* \left(\frac{dl}{(dl, mn)} \right)}{d^{1 + s + v} l^{1 + 2s} m^{1/2 + s} n^{1/2 + v}} \frac{Q'_{dl}(mn)}{(mn)^*{}^2} \frac{\mu(g)}{g^5}.
\end{equation*}
In multiplicative notation we rewrite
\begin{equation}
\label{eq:Zmult}
Z(s, v) = \prod_p \sum_{\substack{(p^g, 2) =1 \\ g \leq 1} } \mathop{ 
\sum_d \sum_l \sum_m \sum_n}_{\substack{m+n \neq 1 \\ d + l + m \leq 1 \\ \text{min}(m+n, g) = 0}}
\frac{(-1)^{d + m } \varphi^* \left(\frac{p^{d+l}}{(p^{d+l}, p^{m+n})} \right)}{p^{d(1 + s +v) + l(1 + 2s) + m(1/2 + s) + n(1/2 + v)}} 
\frac{Q_{p^{d+l}}'(p^{m+n})}{(p^{m+n}){}^*{}^2}\frac{(-1)^g}{p^{5g}}.
\end{equation}
Obviously the Euler product is absolutely convergent for $\text{Re } s > 0$, $\text{Re } v > 0$.  By extracting the lower degree factors we see that
\begin{equation*}
Z(s,v) = \frac{\zeta(1 + 2s)}{\zeta(1 + s + v)} \eta(s,v),
\end{equation*}
where $\eta(s,v)$ is given by an Euler product that is absolutely convergent for $\text{Re } s > -1/6$, $\text{Re } v > -1/6$.  Here we have used Lemma \ref{lem:Q(p^2)} to show that the terms with $m + n = 2$ do not contribute a zero or pole at $s = 0$ or $v = 0$.

Let $I$ be the double integral in the expression for $\Theta$ (with the $2 \pi i$ factor included).  Then
\begin{equation*}
I = \frac{1}{(2 \pi i)^2} \int \int \frac{\zeta(1 + 2s)}{\zeta(1 + s + v)} \frac{M^s}{s^{j +1}} \frac{U^v}{v} \eta(s,v) (2\pi)^{-v} \Gamma(1 + v) G(v) ds dv.
\end{equation*}
To begin we take $\text{Re } s = 1$ and move the line of integration over $v$ to $-\delta$.  The only pole is at $v=0$ and leads to
\begin{equation*}
I = \frac{1}{2 \pi i} \int \frac{\zeta(1 + 2s)}{\zeta(1 + s )} \frac{M^s}{s^{j +1}} \eta(s,0) ds + I_{\delta},
\end{equation*}
where $I_\delta$ is the same integral as $I$ but along $\text{Re } v = -\delta$.  By moving the $s$-line of integration to $\text{Re } s = \varepsilon$ (assume RH for brevity) we obtain
\begin{equation*}
|I_{\delta}| \ll \frac{M^{\varepsilon}}{U^{\delta}}.
\end{equation*}
Since $U$ and $M$ are both fixed powers of $X$ we may take $\varepsilon$ small enough so that $I_{\delta} \ll X^{-\varepsilon}$.  By moving $\text{Re } s < 0$ we obtain an asymptotic formula for $I$.  We clearly obtain
\begin{equation*}
I = \frac{(\log{M})^j}{j!} \onehalf \eta(0,0) + O((\log{M})^{j-1}). %+ X^{-\varepsilon}).
\end{equation*} 
Inserting this formula into $\Theta$, we obtain
\begin{equation*}
\Theta = \frac{AB}{2}  \widehat{w}(0,0) \onehalf \eta(0,0) +O\left(\frac{AB}{\log{X}}\right),
\end{equation*}
using the fact that $\sum_j a_j = P(1) - P(0) = 1$ (by definition of $P(x)$).

All that remains is to compute the value of $\eta(0,0)$.  We use the fact that $\eta(v, v) = Z(v, v)$ and compute with $Z$.  We apply the changes of variable $n \rightarrow n - m$ and $l \rightarrow l - d$ to \eqref{eq:Zmult} and obtain
\begin{equation*}
Z(v, v) = \prod_p \mathop{\sum_{(p^g, 2) =1} \sum_l \sum_{d \leq l}  \sum_{n \neq 1} \sum_{m \leq n} }_{\substack{l + m \leq 1 \\ \text{min}(n, g) = 0}}
\frac{(-1)^{d + m } \varphi^* \left(\frac{p^{l}}{(p^{l}, p^{n})} \right)}{p^{l(1 + 2v) + n(1/2 + v)}} 
\frac{Q_{p^{l}}'(p^{n})}{(p^{n}){}^*{}^2}\frac{(-1)^g}{p^{5g}}.
\end{equation*}
Executing the summation over $d$ first, we obtain zero as long as $d$ can take the value $1$.  Thus we may assume $l = 0$.  Next we sum over $m$ and obtain zero as long as $m$ can take the value $1$.  Hence we may also assume $n = 0$.  
Therefore $Z(v,v) = \zeta^{-1}(5) (1 - 2^{-5})^{-1}$ for $\text{Re } v > 0$.  We conclude that the main term is
\begin{equation*}
\frac14 \zeta^{-1}(5) AB \widehat{w}(0, 0) (1 - 2^{-5})^{-1}.
\end{equation*}

\subsection{The proof of Theorem \ref{thm:LU}}
\label{section:changes}
Since the details of the proof of Theorem \ref{thm:LU} are similar to (but much easier than) the proof of Theorem \ref{thm:LUM}, we simply summarize the differences.  

The analogues of Lemmas \ref{lem:LUMapprox} and \ref{lem:LUMavgclean} are proved in the same way; the absence of the mollifier only simplifies matters. 

The estimation of the remainder term as in Section \ref{section:RT} proceeds essentially unchanged.

Since a different main term arises we shall presently make this calculation.  The main term is the following analogue of \eqref{eq:mollifiedmainterm}
\begin{equation}
\label{eq:mollifiedmainterm2}
\onehalf AB \widehat{w}(0,0)  \sum_{(g, 2) =1} \frac{\mu(g)}{g^5} \sum_{\substack{%(n)_2 \leq X^{\varepsilon} \\ 
(n, g) = 1}} \frac{1}{\sqrt{n}} \frac{Q(n)}{n^*{}^2} Y\left(\frac{2\pi n}{U} \right) \phi\left(\frac{(n)_2}{X^{\varepsilon}}\right).
\end{equation}

We first remove the truncation factor $\phi$ without introducing a new error term.  Let $\Theta'$ be the sum to be computed.  We have
\begin{equation*}
\Theta' = \onehalf AB \widehat{w}(0,0) \frac{1}{(2 \pi i)} \int  Z(v)\frac{U^v}{v} (2\pi)^{-v} \Gamma(1 + v) G(v)  dv,
\end{equation*} 
where
\begin{equation*}
Z(v) =  \sum_{(g, 2) =1} \sum_{(n, g) = 1} \frac{Q(n) }{n^*{}^2 n^{1/2 + v}} \frac{\mu(g)}{g^5}.
\end{equation*}
Using $Q(p^k) = 0$ for $k$ odd (Lemma \ref{lemma:Qpkodd}) we have
\begin{equation*}
Z(v) =  \sum_{(g, 2) =1} \sum_{(n, g) = 1} \frac{Q(n^2) }{n^*{}^2 n^{1 + 2v}} \frac{\mu(g)}{g^5}.
\end{equation*}
Now we simply move the line of integration past $0$ and obtain
\begin{equation*}
\Theta' = \onehalf AB \widehat{w}(0,0) Z(0).
\end{equation*}
We have $Z(0) = \zeta^{-1}(5) (1 - 2^{-5})^{-1} c_S$, where
\begin{equation}
\label{eq:cS}
c_S = \prod_p \left(1 + \left( 1- p^{-5}\right)^{-1}  \sum_{k=1}^{\infty} \frac{Q(p^{2k})}{p^{k + 2}} \right),
\end{equation}
which completes the calculation of the main term.

\section{On the second moments of the partial sums in the approximate functional equation}
\label{section:secondmoments}
There are some difficulties that arise when analyzing the individual partial sums in the approximate functional equation.  This section is devoted to elaborating on this issue.

We prove Propositions \ref{prop:LVsquared}, \ref{prop:LVMsquared}, \ref{prop:Msquared}, and \ref{prop:LV1V2M1M2} in the following sections.  The proofs are similar so to reduce repetition we handle the four cases simultaneously whenever possible.

\subsection{Estimating the remainder terms}
The basic idea of the proof is analogous to that of Theorem \ref{thm:LUM}.  The calculation of the main term becomes significantly more involved.  The restriction to sums of total length less than $5/9$ allows us to much more easily show that the remainder term is small.

We begin by writing
\begin{multline}
\label{eq:LVMEoneHecke}
L_{V_1} M_1(E) L_{V_2} M_2(E) \\
=
 \mathop{\sum_{d_1, d_2} \sum_{l_1, l_2}}_{(d_1d_2l_1l_2, \Delta) = 1} \sum_{m_1, m_2} \sum_{n_1, n_2} \frac{\lambda(m_1n_1)\lambda(m_2n_2) \mu(d_1m_1) \mu(d_2m_2)}{d_1 d_2 l_1 l_2 \sqrt{m_1n_1 m_2 n_2}}
 \mu^2(d_1l_1m_1)\mu^2(d_2l_2m_2)  F,
\end{multline}
where $F$ is shorthand for the following product of test functions
\begin{equation*}
F = Y\left(\frac{2\pi d_1n_1}{V_1} \right) Y\left(\frac{2\pi d_2n_2}{V_2} \right)
P \left( \frac{\log{M_1/d_1m_1l_1^2}}{\log{M_1}} \right)  P \left( \frac{\log{M_2/d_2m_2l_2^2}}{\log{M_2}} \right).
\end{equation*}
With a stretch of the imagination this expression will cover all four desired cases; when $M(E) = 1$ we simply set $d_i = l_i = m_i = 1$ and set $P(x) = 1$.  Likewise, when $L_V = 1$ we set $d_i = l_i = n_i = 1$ and set $Y(u) = 1$.

Applying the Hecke relation again we obtain
\begin{multline}
\label{eq:LVME}
L_{V_1} M_1(E) L_{V_2} M_2(E) \\
=
\mathop{\sum_{d_1, \cdots, n_2} \sum_{f | (m_1n_1, m_2 n_2)}}_{(d_1d_2l_1l_2f, \Delta) = 1} \frac{\lambda(m_1n_1m_2n_2 f^{-2})}{d_1 d_2 l_1 l_2 \sqrt{m_1n_1 m_2 n_2}}
\mu(d_1m_1) \mu(d_2m_2) \mu^2(d_1l_1m_1)\mu^2(d_2l_2m_2) F.
\end{multline}

Now we sum over the family.  The only factor depending on $E$ in the above expression is $\lambda(m_1n_1m_2n_2 f^{-2}) \psi_{\Delta}(d_1l_1d_2 l_2f)$.  Let
\begin{equation*}
T = \sum_{(a, b) \in S} \lambda_{a, b}(m_1n_1m_2n_2 f^{-2}) \psi_{\Delta}(d_1l_1d_2 l_2f) w_X(a,b).
\end{equation*}
Then
\begin{equation*}
T = \sum_{c | d_1 d_2 l_1 l_2} \mu(c) \sum_{\substack{(a, b) \in S \\ \Delta \equiv 0 \shortmod{c}}} \lambda_{a, b}(m_1n_1m_2n_2 f^{-2})  w_X(a,b).
\end{equation*}
By another application of M\"{o}bius inversion, we have
\begin{equation*}
T = \sum_{c | d_1 d_2 l_1 l_2} \mu(c) \sum_{(g, 2m_1m_2n_1n_2f^{-2}) = 1} \mu(g) \sum_{\substack{a, b \\ g^6 \Delta \equiv 0 \shortmod{c}}} \lambda_{ag^2, bg^3}(m_1n_1m_2n_2 f^{-2})  w_X(ag^2,bg^3).
\end{equation*}
To the inner sum we apply \eqref{eq:usefulZ} with $r = m_1 n_1 m_2 n_2 f^{-2}$.  The zero frequencies $h = k = 0$ lead to a main term, say $T_0$.  The terms with $h \neq 0$ or $k \neq 0$ are estimated trivially.  Writing $T = T_0 + R$, we have the trivial bound
\begin{equation*}
R \ll AB X^{\varepsilon} r_2^*{}^2 \sum_{c | d_1 d_2 l_1 l_2} \sum_{g \ll X^{1/6}} \frac{r_0 \sqrt{c_0}}{q^2 g^5} \left( \frac{g^2q}{A} + \frac{g^3 q}{B} + \frac{g^5 q^2}{AB} \right),
\end{equation*}
using Weil's bound for the complete character/exponential sums appearing in \eqref{eq:usefulZ}.  Recall that $r = r_1 r_2$, where $r_1$ is the product of primes exactly dividing $r$, $c = (c, g) c_0 c_1 c_2$, $c_1 | r_1$, $c_2 | r_2$, $(r, c_0) = 1$, and $q = r_1 r_2^* c_0$.  All we need is that $c_0 \leq c$ and that $q \leq r_1 r_2^* c$ to see that
\begin{equation*}
R \ll AB r_1 r_2^*{}^2 X^{\varepsilon} \left(\frac{1}{A r_1 r_2^*} + \frac{1}{B r_1 r_2^*} + \frac{\sqrt{d_1 d_2 l_1 l_2} X^{1/6}}{AB}\right).
\end{equation*}
Inserting this bound for $R$ into the expression for the average value of $L_V M(E) L_{V'} M'(E)$ shows that this remainder term contributes $O(AB X^{-\delta})$ for $V_1V_2M_1M_2 \ll X^{5/9 - \varepsilon}$.

We now concentrate on $T_0$.  The main term is given by \eqref{eq:mainterm}, hence
\begin{equation*}
T_0 = AB \widehat{w}(0,0) \sum_{c | d_1 d_2 l_1 l_2} \mu(c) \sum_{(g, 2m_1m_2n_1n_2f^{-2}) = 1} \frac{\mu(g)}{g^5} \delta_{\square}(m_1 m_2 n_1 n_2 f^{-2}) \frac{Q_c(m_1m_2n_1n_2f^{-2})}{r^*{}^2},
\end{equation*}
where $\delta_{\square}$ is the characteristic function of squares.  Reversing the M\"{o}bius inversion in $c$, we obtain
\begin{equation}
\label{eq:T0}
T_0 = AB \widehat{w}(0,0) \sum_{(g, 2m_1m_2n_1n_2f^{-2}) = 1} \frac{\mu(g)}{g^5} \delta_{\square}(m_1 m_2 n_1 n_2 f^{-2}) \frac{Q_{d_1d_2l_1l_2}'(m_1m_2n_1n_2f^{-2})}{r^*{}^2}.
\end{equation}
Now that we have the main term in this form we can use zeta-function methods (as in Section \ref{section:mainterm}) to evaluate the summation over the various variables present.

\subsection{Evaluating the main terms}

The main term to be evaluated is simply \eqref{eq:T0} inserted into \eqref{eq:LVME}.  Precisely, it is
\begin{multline*}
AB \widehat{w}(0,0) \sum_{d_1, d_2} \sum_{l_1, l_2} \mathop{\sum_{m_1, m_2} \sum_{n_1, n_2}}_{m_1 m_2 n_1 n_2 = \square} \sum_{f|(m_1 n_1, m_2 n_2)} \sum_{(g, 2m_1m_2n_1n_2f^{-2}) = 1} \frac{\mu(g)}{g^5} \\
\cdot
\frac{\mu(d_1m_1) \mu(d_2m_2) \mu^2(d_1l_1m_1)\mu^2(d_2l_2m_2)}{d_1 d_2 l_1 l_2 \sqrt{m_1n_1 m_2 n_2}}    \frac{Q_{d_1d_2l_1l_2}'(m_1m_2n_1n_2f^{-2})}{(m_1 m_2 n_1 n_2 f^{-2})^*{}^2} F.
\end{multline*}
We apply Mellin inversion on $F$ and write the sum as an integral as follows
\begin{multline}
\label{eq:maintermintegral}
AB \widehat{w}(0,0) \sum_{j_1} \sum_{j_2} \frac{a_{j_1} a_{j_2} j_1! j_2!}{(\log{M_1})^{j_1}(\log{M_2})^{j_2}} \frac{1}{(2 \pi i)^4} 
\\
\cdot
\int \int \int \int 
Z(s_1, s_2, v_1, v_2) 
\frac{M_1^{s_1}}{s_1^{j_1 + 1}}
\frac{M_2^{s_2}}{s_2^{j_2 + 1}} 
\frac{V_1^{v_1}}{v_1}
\frac{V_2^{v_2}}{v_2}
 \frac{G(v_1)}{(2\pi)^{v_1}} \frac{G(v_2)}{(2\pi)^{v_2}} \Gamma(v_1 + 1) \Gamma(v_2 + 1) ds_1 ds_2 dv_1 dv_2,
\end{multline}
where
\begin{multline*}
Z(s_1, s_2, v_1, v_2) = \sum_{d_1, d_2} \sum_{l_1, l_2} \mathop{\sum_{m_1, m_2} \sum_{n_1, n_2}}_{m_1 m_2 n_1 n_2 = \square} \sum_{f|(m_1 n_1, m_2 n_2)} \sum_{(g, 2m_1m_2n_1n_2f^{-2}) = 1} \frac{\mu(g)}{g^5} 
\\
\cdot
\frac{\mu(d_1m_1) \mu(d_2m_2) \mu^2(d_1l_1m_1)\mu^2(d_2l_2m_2)}{d_1^{1 + s_1 + v_1} d_2^{1 + s_2 + v_2} l_1^{1 + 2s_1} l_2^{1 + 2s_2} m_1^{1/2 + s_1} m_2^{1/2 + s_2} n_1^{1/2 + v_1} n_2^{1/2 +v_2 }}    \frac{Q_{d_1d_2l_1l_2}'(m_1m_2n_1n_2f^{-2})}{(m_1 m_2 n_1 n_2 f^{-2})^*{}^2}.
\end{multline*}
In multiplicative form we have
\begin{multline*}
Z = \prod_{p \neq 2}  \mathop{\sum_{d_1, d_2} \sum_{l_1, l_2} \sum_{m_1, m_2} \sum_{n_1, n_2}}_{\substack{m_1 + m_2+ n_1+ n_2 \equiv 0 \shortmod{2} \\ d_1 + l_1 + m_1 \leq 1 \\ d_2 + l_2 + m_2 \leq 1}} \mathop{\sum \sum}_{\substack{f \leq \text{min}(m_1+ n_1, m_2 + n_2) \\ \text{min}(g, m_1+m_2+n_1+ n_2 -2 f) = 0}} \frac{\mu(p^g)}{p^{5g}} 
\\
\cdot
\frac{(-1)^{d_1 + d_2 + m_1 + m_2} Q_{p^{d_1+d_2+l_1+l_2}}'(p^{m_1+m_2+n_1+n_2 -2 f}) [(p^{m_1+ m_2 + n_1 + n_2 -2 f})^*{}]^{-2}}{p^{d_1(1 + s_1 + v_1) + d_2(1 + s_2 + v_2) + l_1(1 + 2s_1)+ l_2(1 + 2s_2) + m_1(1/2 + s_1) + m_2(1/2 + s_2) + n_1(1/2 + v_1) +  n_2(1/2 +v_2)}}.
\end{multline*}
In order to understand the behavior of $Z$ near $(0, 0, 0, 0)$ we extract the lower-degree factors of $Z$.  The only factors in the Euler product that have an effect near $0$ are the coefficients of $p^{-1}$, so to speak.  It is at this point that our treatments of Propositions 
\ref{prop:LVsquared}, \ref{prop:LVMsquared}, \ref{prop:Msquared}, and \ref{prop:LV1V2M1M2} diverge.  In each case we shall reduce the computation of the main term to the asymptotic evaluation of a certain multidimensional integral.  The computations of the necessary integrals are performed in Section \ref{section:integrals}.

{\bf The case $M(E) = 1$, $V_1 = V_2 = V$}.  The analysis of this case will prove Proposition \ref{prop:LVsquared}.  The expression \eqref{eq:maintermintegral} simplifies to the form
\begin{equation*}
AB \widehat{w}(0,0)  \frac{1}{(2 \pi i)^2} 
\int \int 
Z(v_1, v_2) 
 \frac{V^{v_1 + v_2}}{v_1 v_2}
 \frac{G(v_1)}{(2\pi)^{v_1}} \frac{G(v_2)}{(2\pi)^{v_2}} \Gamma(v_1 + 1) \Gamma(v_2 + 1) dv_1 dv_2,
\end{equation*}
where
\begin{equation*}
Z(v_1, v_2) = \prod_{p \neq 2} \mathop{\sum \sum}_{\substack{n_1, n_2 \\ n_1+ n_2 \equiv 0 \shortmod{2}}} \mathop{\sum \sum}_{\substack{f \leq \text{min}(n_1, n_2) \\ \text{min}(g, n_1+ n_2 -2 f) = 0}} \frac{\mu(p^g)}{p^{5g}}  
\frac{ Q(p^{n_1+n_2 -2 f}) [(p^{n_1 + n_2 -2 f})^*{}]^{-2}}{p^{n_1(1/2 + v_1) +  n_2(1/2 +v_2)}}.
\end{equation*}
The case $n_1 = n_2 = 1$, $f = 1$ leads to the only term of low degree.  We therefore have
\begin{equation*}
Z(v_1, v_2) = \zeta(1 + v_1 + v_2) \eta_1(v_1, v_2),
\end{equation*}
where $\eta_1$ is an Euler product absolutely convergent in some region $\text{Re } v_j > -\delta$ for $\delta > 0$.  %Now we evaluate the integral by moving the lines of integration.  
The desired main term is then given by
\begin{equation*}
AB \widehat{w}(0,0)  \frac{1}{(2 \pi i)^2} 
\int \int 
\zeta(1 + v_1 + v_2)  
\frac{V^{v_1 + v_2}}{v_1 v_2}
\eta_1(v_1, v_2) \frac{G(v_1)}{(2\pi)^{v_1}} \frac{G(v_2)}{(2\pi)^{v_2}} \Gamma(v_1 + 1) \Gamma(v_2 + 1) dv_1 dv_2.
\end{equation*}
This integral is computed with Proposition \ref{prop:integralone}; we obtain
\begin{equation*}
\eta_1(0, 0) \log{V} AB \widehat{w}(0,0) + O(AB).
\end{equation*}
Here $\eta_1(0,0)$ is an arithmetical factor that can be given by an Euler product similarly to \eqref{eq:cS} if desired.  
This concludes the proof of Proposition \ref{prop:LVsquared}.

{\bf The case $L_V = 1$, $M_1 = M_2 = M$}.  Here we shall prove Proposition \ref{prop:Msquared}.  In this case the expression \eqref{eq:maintermintegral} takes the form
\begin{equation*}
AB \widehat{w}(0,0) \sum_{j_1} \sum_{j_2} \frac{a_{j_1} a_{j_2} j_1! j_2!}{(\log{M})^{j_1+ j_2}} \frac{1}{(2 \pi i)^2} 
 \int \int 
Z(s_1, s_2) 
\frac{M^{s_1 + s_2}}{s_1^{j_1 + 1}s_2^{j_2 + 1}}
ds_1 ds_2,
\end{equation*}
where
%\begin{equation*}
%Z(s_1, s_2) = \prod_{p}  \mathop{ \sum_{m_1, m_2} }_{\substack{m_1 + m_2 \equiv 0 \shortmod{2}}} \mathop{\sum \sum}_{\substack{f \leq %\text{min}(m_1, m_2) \\ \text{min}(g, m_1+m_2 -2 f) = 0}} \frac{\mu(g)}{g^5} 
%\frac{(-1)^{m_1 + m_2} Q(p^{m_1+m_2 -2 f})}{(p^{m_1+ m_2  -2 f})^*{}^{2} p^{m_1(1/2 + s_1) + m_2(1/2 + s_2) }}.
%\end{equation*}
\begin{multline*}
Z(s_1, s_2) = \prod_{p \neq 2}  \mathop{\sum_{l_1, l_2} \sum_{m_1, m_2}}_{\substack{m_1 + m_2 \equiv 0 \shortmod{2} \\ l_1 + m_1 \leq 1 \\ l_2 + m_2 \leq 1}} \mathop{\sum \sum}_{\substack{f \leq \text{min}(m_1, m_2 ) \\ \text{min}(g, m_1+m_2 -2 f) = 0}} \frac{\mu(p^g)}{p^{5g}}  
\\
\cdot
\frac{(-1)^{ m_1 + m_2} Q_{p^{l_1+l_2}}'(p^{m_1+m_2-2 f}) [(p^{m_1+ m_2 -2 f})^*{}]^{-2}}{p^{l_1(1 + 2s_1)+ l_2(1 + 2s_2) + m_1(1/2 + s_1) + m_2(1/2 + s_2)}}.
\end{multline*}

There are three terms that lead to a low degree in the Euler product.  The cases are $l_1 = 1$, $l_2 = m_1 = m_2 =0$; $l_2 = 1$, $l_1 = m_1 = m_2 = 0$; and $m_1 = m_2 = 1$, $f=1$, $l_1 = l_2 = 0$. Hence
\begin{equation*}
Z(s_1, s_2) = \zeta(1 + s_1 + s_2) \zeta(1 + 2s_1) \zeta(1 + 2s_2) \eta_2(s_1, s_2),
\end{equation*}
where $\eta_2$ has the same properties as $\eta_1$.  The integral we need to compute is
\begin{equation}
\label{eq:Iintegral}
%I = 
\frac{1}{(2 \pi i)^2} 
 \int \int 
\zeta(1 + s_1 + s_2) \zeta(1 + 2s_1) \zeta(1 + 2s_2)
\frac{M^{s_1 + s_2}}{s_1^{j_1 + 1}s_2^{j_2 + 1}} \eta_2(s_1, s_2)
ds_1 ds_2,
\end{equation}
which is accomplished with Proposition \ref{prop:integraltwo}.

Using the computation and summing over $j_1$ and $j_2$ leads to the main term being
\begin{equation*}
\frac14 \eta_2(0, 0) (\log{M})^3 AB \widehat{w}(0, 0) \sum_{j_1} \sum_{j_2} \frac{(j_1 + 1) (j_2 + 1) a_{j_1} a_{j_2}}{j_1 + j_2 + 3}.
\end{equation*}
Setting $F(x) = x \frac{d}{dx} xP(x) = x^2P'(x) + xP(x)$ we have
\begin{equation*}
\sum_{j_1} \sum_{j_2} \frac{(j_1 + 1) (j_2 + 1) a_{j_1} a_{j_2}}{j_1 + j_2 + 3} \int_0^1 F(x)^2 dx.
\end{equation*}

Hence the main term is
\begin{equation*}
\frac14 \eta_2(0, 0) (\log{M})^3 AB \widehat{w}(0,0) \int_0^1 F(x)^2 dx.
\end{equation*}
The value $\eta_2(0,0)$ is an absolute arithmetical constant that can be given by a formula similar to \eqref{eq:cS}.  This concludes the proof of Proposition \ref{prop:Msquared}.

{\bf The general case}.  We simultaneously handle the cases relevant for Propositions \ref{prop:LVMsquared} and \ref{prop:LV1V2M1M2} for a while longer before treating the cases differently.  

In the analysis of $Z(s_1, s_2, v_1, v_2)$, there are eight ways to get a low degree power of $p$: 1) $d_1 = 1$ with the other variables zero (we omit mention of this condition in the remaining cases although of course it holds)); 2) $d_2 = 1$; 3) $l_1 = 1$; 4) $l_2 = 1$; 5) $m_1  = m_2 = 1$, $n_1 = n_2 = 0$; 6) $m_1 = n_2 = 1$, $n_1 = m_2 = 0$; 7) $n_1 = m_2 = 1$, $m_1  = n_2 = 0$ ; and 8) $m_1 = m_2 = 0$, $n_1 = n_2 = 1$.  We therefore have
\begin{equation*}
Z(s_1, s_2, v_1, v_2) = \frac{\zeta(1 + 2s_1)\zeta(1 + 2s_2)\zeta(1 + s_1 + s_2)\zeta(1 + v_1 + v_2)}{\zeta(1 + s_1 + v_1)\zeta(1 + s_2 + v_2)\zeta(1 + s_1 + v_2)\zeta(1 + s_2 + v_1)} \eta(s_1, s_2, v_1, v_2),
\end{equation*}
where $\eta$ is an Euler product absolutely convergent in some region $\text{Re } s_i, v_j > -\delta$ for $\delta > 0$.
The main term takes the form
\begin{equation*}
AB \widehat{w}(0,0) \sum_{j_1} \sum_{j_2} \frac{a_{j_1} a_{j_2} j_1! j_2!}{(\log{M_1})^{j_1}(\log{M_2})^{j_2}} I,
\end{equation*}
where
\begin{multline*}
I = \frac{1}{(2 \pi i)^4} \int \int \int \int \frac{\zeta(1 + 2s_1)\zeta(1 + 2s_2)\zeta(1 + s_1 + s_2)\zeta(1 + v_1 + v_2)}{\zeta(1 + s_1 + v_1)\zeta(1 + s_1 + v_2)\zeta(1 + s_2 + v_1)\zeta(1 + s_2 + v_2)} 
\\
\cdot \frac{M_1^{s_1}}{s_1^{j_1 + 1}}
\frac{M_2^{s_2}}{s_2^{j_2 + 1}} 
\frac{V_1^{v_1}}{v_1}
\frac{V_2^{v_2}}{v_2}
\eta(s_1, s_2, v_1, v_2)
 \frac{G(v_1)}{(2\pi)^{v_1}} \frac{G(v_2)}{(2\pi)^{v_2}} \Gamma(v_1 + 1) \Gamma(v_2 + 1) ds_1 ds_2 dv_1 dv_2.
\end{multline*}
This integral is asymptotically computed with Proposition \ref{prop:integralthree}; there is different behavior depending on if $V_1 =V_2$ or $V_1 \neq V_2$.  We easily obtain the results of Propositions \ref{prop:LVMsquared} and \ref{prop:LV1V2M1M2}
by applying Proposition \ref{prop:integralthree}. %and summing over $j_1$ and $j_2$.

%Gathering our main terms, we obtain that the asymptotic in Proposition \ref{prop:LV1V2M1M2} is
%\begin{equation*}
%AB \widehat{w}(0,0) c_0 c(M_1, M_2) + O\left(\frac{AB}{\log{X}}\right),
%\end{equation*}
%where $c_0$ is an absolute arithmetical constant and $c(M_1, M_2)$ is a combinatorial constant depending on the lengths of the mollifiers as well as the particular polynomial $P(x)$ given by \eqref{eq:mollifier}.  In special cases such as if $M_1 = M_2$ or if $P(x) = x$ it is possible to simplify $c(M_1, M_2)$.  This completes the proof of Proposition \ref{prop:LV1V2M1M2}.

\subsection{Integral computations}
\label{section:integrals}

\begin{myprop}
\label{prop:integralone}
Let $g(v_1, v_2)$ be a holomorphic function in $\text{Re } v_1$, $\text{Re } v_2 \geq -\delta$, for some $\delta > 0$, and suppose $g(v_1, v_2)  \ll (1 + |v_1|)^{-2} (1 + |v_2|)^{-2}$ for $-\delta \leq \text{Re } v_i \leq 1$, $i=1, 2$.  Let
%$v_1^{\alpha} v_2^{\beta}g^{(\alpha, \beta)}(v_1, v_2)  \ll_{\alpha, \beta} (1 + |v_1|)^{-2} (1 + |v_2|)^{-2}$, the superscript indicating partial differentiation.  Let
\begin{equation*}
I = \frac{1}{(2 \pi i)^2} 
\int_{(1)} \int_{(1)}
\zeta(1 + v_1 + v_2)  
\frac{V^{v_1 + v_2}}{v_1 v_2}
g(v_1, v_2) dv_1 dv_2.
\end{equation*}
Then
\begin{equation*}
I = g(0,0) \log{V} + O(1)
\end{equation*}
as $V \rightarrow \infty$.
\end{myprop}

\begin{proof}
We begin by taking the lines of integration at $\text{Re } v_1 =  \text{Re } v_2 = \varepsilon$ with $0 < \varepsilon < \delta/3$.  Moving the line of integration over $v_2$ to $\text{Re } v_2 = - 2\varepsilon$ we pick up poles at $v_2 = 0$ and at $v_2 = -v_1$. The contribution from the new line of integration is $\ll V^{-\varepsilon}$.  

The $v_2 = -v_1$ pole gives
\begin{equation*}
\frac{1}{2 \pi i} 
\int_{(\varepsilon)}  
\frac{-1}{v_1^2}
g(v_1, -v_1)  dv_1 \ll 1.
\end{equation*}

The pole at $v_2 = 0$ leads to
\begin{equation*}
 \frac{1}{2 \pi i} 
\int_{(\varepsilon)}  
\zeta(1 + v_1)  
 \frac{V^{v_1}}{v_1}
g(v_1, 0)  dv_1 .
\end{equation*}
Now we move $v_1$ to the left of the imaginary axis.  A power of $V$ is saved in the integration over the new line. There is a double pole at $v_1 = 0$, so the contribution from this pole is $g(0,0) \log{V} + O(1)$.
\end{proof}

\begin{myprop}
\label{prop:integraltwo}
Let $g$ be as in Proposition \ref{prop:integralone} and set
\begin{equation*}
I = \frac{1}{(2 \pi i)^2} 
 \int_{(1)} \int_{(1)} 
\zeta(1 + s_1 + s_2) \zeta(1 + 2s_1) \zeta(1 + 2s_2)
\frac{M^{s_1 + s_2}}{s_1^{j_1 + 1}s_2^{j_2 + 1}} g(s_1, s_2)
ds_1 ds_2.
\end{equation*}
Then
\begin{equation*}
I = \frac14 g(0, 0)  \frac{(\log{M})^{j_1 + j_2 + 3} (j_1 + 1)! (j_2 + 1)!}{j_1 + j_2 + 3} \left(1 + O\left(\frac{1}{\log{M}}\right) \right)
\end{equation*}
\end{myprop}
\begin{proof}
It simplifies the calculation of the integral to separate the variables $s_1$ and $s_2$ by means of the Dirichlet series expansion for $\zeta(1 + s_1 + s_2)$.  We may freely interchange the summation and integration in the region of absolute convergence $\text{Re } s_1 > 0$, $\text{Re } s_2 > 0$, so we obtain
\begin{equation*}
I =  \sum_{n } \frac{1}{n} \frac{1}{(2 \pi i)^2} 
 \int_{(1)} \int_{(1)}  
\frac{\zeta(1 + 2s_1)}{s_1^{j_1 + 1}} \frac{\zeta(1 + 2s_2)}{s_2^{j_2 + 1}}(M/n)^{s_1 + s_2} g(s_1, s_2)
ds_1 ds_2.
\end{equation*}
It is obvious that the terms with $n > M$ contribute $O(1)$ to $I$.  For $n \leq M$ we compute the integral asymptotically by moving the lines of integration past the origin, obtaining
\begin{equation*}
I = \frac14 g(0, 0) (j_1 + 1)! (j_2 + 1)!\sum_{n \leq M} \frac{(\log{(M/n)})^{j_1 + j_2 + 2}}{n} \left(1 + O\left(\frac{1}{\log{M}}\right) \right) + O(1).
\end{equation*}
An application of Euler-Maclaurin shows
\begin{align*}
I & = \frac14 g(0, 0) (j_1 + 1)! (j_2 + 1)! \int_{1}^{M}  \frac{(\log{(M/t)})^{j_1 + j_2 + 2}}{t} dt \left(1 + O\left(\frac{1}{\log{M}}\right) \right) + O(1)
\\
& = \frac14 g(0, 0)  \frac{(\log{M})^{j_1 + j_2 + 3} (j_1 + 1)! (j_2 + 1)!}{j_1 + j_2 + 3} \left(1 + O\left(\frac{1}{\log{M}}\right) \right).
\end{align*}
\end{proof}

\begin{myprop}
\label{prop:integralthree}
Let $g(s_1, s_2, v_1, v_2)$ be holomorphic in $\text{Re } s_i, v_i \geq -\delta$, for some $\delta > 0$, and suppose 
$g(s_1, s_2, v_1, v_2) \ll (1 + |t_1|)^{-2}(1 + |t_2|)^{-2}(1 + |v_1|)^{-2}(1 + |v_1|)^{-2}$ provided $-\delta \leq \text{Re } v_i \leq 1$ and $\text{Re } s_i \geq -\delta$, $i=1,2$, and where $t_i = \text{Im } s_i$.  Furthermore suppose that $g(0,0,-it,it) > 0$ for all $t \in \mr$.  
Let $V_1, V_2, M_1$, and $M_2$ all be fixed powers of a parameter $X$ and set
\begin{multline*}
I = \frac{1}{(2 \pi i)^4} \int_{(1)} \int_{(1)} \int_{(1)}\int_{(1)} \frac{\zeta(1 + 2s_1)\zeta(1 + 2s_2)\zeta(1 + s_1 + s_2)\zeta(1 + v_1 + v_2)}{\zeta(1 + s_1 + v_1)\zeta(1 + s_1 + v_2)\zeta(1 + s_2 + v_1)\zeta(1 + s_2 + v_2)} 
\\
\cdot \frac{M_1^{s_1}}{s_1^{j_1 + 1}}
\frac{M_2^{s_2}}{s_2^{j_2 + 1}} 
\frac{V_1^{v_1}}{v_1}
\frac{V_2^{v_2}}{v_2}
g(s_1, s_2, v_1, v_2)
 ds_1 ds_2 dv_1 dv_2.
\end{multline*}
If $V_1 = V_2$ then there exists a positive constant $c$ such that
\begin{equation*}
I \sim c (\log{X})^{j_1 + j_2 +3}
\end{equation*}
as $X \rightarrow \infty$, 
and if $V_1 \neq V_2$ then there exists a positive constant $c'$ such that
\begin{equation*}
I \sim c' (\log{X})^{j_1 + j_2} 
\end{equation*}
as $X \rightarrow \infty$.
\end{myprop}
The proof gives explicit computations of $c$ and $c'$.
\begin{proof}
This material is very similar to that of \cite{KMVdK}, Section 5.  We begin with all the lines of integration very close but to the right of the imaginary axes.  Without loss of generality we may assume $V_1 \leq V_2$ and $M_1 \leq M_2$.  We first move the line of integration over $v_2$ to a contour $\gamma_2$ of the form $-C (\log{(|t| + 2)})^{-1} + i t$, for some $C > 0$, thereby picking up a pole only at $v_2 = 0$ (assuming $C > 0$ is small enough).  In what follows we say that a contour is of the type $\gamma$ if it has the same form as $\gamma_2$ (but with perhaps a different constant).  We obtain
\begin{equation*}
I = I_1 + I_{\gamma},
\end{equation*}
where
\begin{multline*}
I_1 = \frac{1}{(2 \pi i)^3}  \int \int \int \frac{\zeta(1 + s_1 + s_2)\zeta(1 + v_1)}{\zeta(1 + s_1 + v_1)\zeta(1 + s_2 + v_1)} 
\frac{M_1^{s_1}}{s_1^{j_1 + 1}}
\frac{M_2^{s_2}}{s_2^{j_2 + 1}} 
\frac{V_1^{v_1}}{v_1}
g_1(s_1, s_2, v_1)
   ds_1 ds_2 dv_1, %\\
 %+ I_{\gamma},
\end{multline*}
\begin{equation*}
g_1(s_1, s_2, v_1) = \frac{\zeta(1 + 2s_1) \zeta(1 + 2s_2)}{\zeta(1 + s_1) \zeta(1 + s_2)} g(s_1, s_2, v_1, 0),
\end{equation*}
and $I_\gamma$ is given by the same formula as $I$ but with the contour of integration over $v_2$ along $\gamma_2$.  %Let $I_1$ be the triple integral above.  
We defer the treatment of $I_\gamma$ and continue with $I_1$.  

First we move the line of integration over $v_1$ to a contour $\gamma_1$ of type $\gamma$.  We need to extract the coefficient of $v_1^{-1}$ in the Laurent series expansion of the integrand.  Set $h(s) = [\frac{1}{\zeta}]'(1+s)$ (note $h(s)= 1 + O(s)$ for $s$ small).  The relevant expansions are 
\begin{align*}
\zeta(1 + v_1) & = v_1^{-1} +c + ... \\
(\zeta(1 + s + v_1))^{-1} & = (\zeta(1  + s))^{-1} + v_1 h(s) + ... \\ 
V_1^{v_1} & = 1 + v_1 \log{V_1} + ....
\end{align*}
Hence we obtain
\begin{multline}
\label{eq:I2I3}
I_1 \sim \log{V_1}\frac{1}{(2 \pi i)^2} \int \int \frac{\zeta(1 + s_1 + s_2)}{\zeta(1 + s_1)\zeta(1 + s_2)} 
\frac{M_1^{s_1}}{s_1^{j_1 + 1}}
\frac{M_2^{s_2}}{s_2^{j_2 + 1}} 
g_1(s_1, s_2, 0)
  ds_1 ds_2 \\
+  \frac{1}{(2 \pi i)^2} \int \int \zeta(1 + s_1 + s_2) \left( \frac{h(s_1)}{\zeta(1 + s_2)} + \frac{h(s_2)}{\zeta(1 + s_1)} \right)
\frac{M_1^{s_1}}{s_1^{j_1 + 1}}
\frac{M_2^{s_2}}{s_2^{j_2 + 1}} 
g_1(s_1, s_2, 0)
  ds_1 ds_2 + I_{1, \gamma_1},
\end{multline}
where $I_{1, \gamma_1}$ is the same triple integral as $I_1$ but with $v_1$ integrated along $\gamma_1$.   

We claim that $I_{1, \gamma_1}$ is small (with an arbitrary power of $\log{X}$ saving).  To see this, move the integrations over $s_1$ and $s_2$ to contours of the type $\gamma_r$, where $\gamma_r$ denotes a contour of the same type as $\gamma$ except reflected through the imaginary axis.  Here we take the contours sufficiently close to the imaginary axis so that the savings in $V_1$ win over the contribution from $M_1$ and $M_2$.  Here we simply use the classical bound  $\zeta^{-1}(1 + s) \ll \log(2 + |t|)$ along a contour of the type $\gamma$ with $C$ sufficiently small (see (3.11.8) of \cite{Titchmarsh}). The situation where $s_1 + s_2 \approx 0$ is easily handled by trivial estimations.

Set $I_1 \sim \log{V_1} \; I_2 + I_3$ in the way indicated by \eqref{eq:I2I3}.  We compute $I_2$ first.  Using the Dirichlet series expansion of $\zeta(1 + s_1 + s_2)$, we have
\begin{equation*}
I_2 = \sum_{n} \frac{1}{n} \frac{1}{(2 \pi i)^2} \int \int \frac{1}{\zeta(1 + s_1)\zeta(1 + s_2)} 
\frac{(M_1/n)^{s_1}}{s_1^{j_1 + 1}}
\frac{(M_2/n)^{s_2}}{s_2^{j_2 + 1}} 
g_1(s_1, s_2, 0)
  ds_1 ds_2.
\end{equation*}
The integral vanishes unless $n \leq \text{min}(M_1, M_2) = M_1$, which can be seen by taking the line of integration over $s_1$ arbitrarily far to the right (here we use the fact that $g$ is bounded in terms of $s_1$ in the right half plane).  Thus we obtain
\begin{equation*}
I_2 \sim \frac{g(0)}{(j_1-1)! (j_2 - 1)!}\sum_{n \leq M_1} \frac{1}{n} (\log({M_1/n}))^{j_1 - 1} (\log({M_2/n}))^{j_2 - 1},
\end{equation*}
by moving the paths of integration to contours of the type $\gamma$. 

An application of Euler-Maclaurin shows
\begin{align*}
I_2 &\sim \frac{g(0)}{(j_1-1)! (j_2 - 1)!} \int_{1}^{M_1} \frac{1}{t} (\log({M_1/t}))^{j_1 - 1} (\log({M_2/t}))^{j_2 - 1} dt \\
&= \frac{g(0)}{(j_1-1)! (j_2 - 1)!} \int_{1}^{M_1} (\log{t})^{j_1 - 1} (\log{(tM_2/M_1)})^{j_2 - 1} \frac{dt}{t} \\
& = \frac{g(0)}{(j_1-1)! (j_2 - 1)!} \sum_{k = 0}^{j_2 - 1} \binom{j_2 -1}{k} (\log{(M_2/M_1)})^{k} \int_{1}^{M_1} (\log{t})^{j_1 +  j_2 - k- 2} \frac{dt}{t} \\
& = \frac{g(0)}{(j_1-1)! (j_2 - 1)!} \sum_{k = 0}^{j_2 - 1} \binom{j_2 -1}{k}  \frac{(\log{(M_2/M_1)})^{k}(\log{M_1})^{j_1 + j_2 -k -1} }{j_1 + j_2 - k - 1}.
\end{align*}
Using $M_1 = X^{\beta_1}$ and $M_2 = X^{\beta_2}$, this simplifies as
\begin{equation}
\label{eq:I2}
I_2 \sim \frac{g(0) (\log{X})^{j_1 + j_2 - 1}}{(j_1-1)! (j_2 - 1)!} \sum_{k = 0}^{j_2 - 1} \binom{j_2 -1}{k}
\frac{(\beta_2 - \beta_1)^k \beta_1^{j_1 + j_2 - k - 1} }{j_1 + j_2 - k - 1}.
\end{equation}

We now compute $I_3$.  Separating the variables as in the calculation of $I_2$, we have
\begin{align*}
I_3 &= \sum_{n \leq M_1} \frac{1}{n} \frac{1}{(2 \pi i)^2} \int \int \left( \frac{h(s_1)}{\zeta(1 + s_2)} + \frac{h(s_2)}{\zeta(1 + s_1)} \right)\frac{(M_1/n)^{s_1}}{s_1^{j_1 + 1}}
\frac{(M_2/n)^{s_2}}{s_2^{j_2 + 1}} 
g_1(s_1, s_2, 0)
  ds_1 ds_2 \\
& \sim \frac{g(0)}{j_1! j_2!} \sum_{n \leq M_1} \frac{1}{n} \left( j_2(\log{M_1/n})^{j_1} (\log({M_2/n}))^{j_2 - 1} + j_1(\log({M_1/n}))^{j_1 -1} (\log({M_2/n}))^{j_2}\right) \\
& \sim \frac{g(0)}{j_1! j_2!} \int_{1}^{M_1} \left(j_2 (\log{t})^{j_1} (\log{(tM_2/M_1)})^{j_2 - 1}  + j_1(\log{t})^{j_1 - 1} (\log{(tM_2/M_1)})^{j_2}\right)\frac{dt}{t}.
\end{align*}
We have already computed this integral while computing $I_2$, so we immediately have
\begin{multline}
\label{eq:I3}
I_3 \sim \frac{g(0) (\log{X})^{j_1 + j_2}}{j_1! j_2!} \left(j_2 \sum_{k = 0}^{j_2 - 1} \binom{j_2 -1}{k}
\frac{(\beta_2 - \beta_1)^k \beta_1^{j_1 + j_2 - k } }{j_1 + j_2 - k } \right. \\
\left. + j_1 \sum_{k = 0}^{j_2 } \binom{j_2 }{k}
\frac{(\beta_2 - \beta_1)^k \beta_1^{j_1 + j_2 - k } }{j_1 + j_2 - k }\right).
\end{multline}
We have computed these integrals to higher precision than we need because these computations may be of use.  Notice that when $\beta_1 = \beta_2$ the formulas simplify as
\begin{align*}
I_2 & \sim \frac{g(0) (\log{M})^{j_1 + j_2 - 1}}{(j_1-1)! (j_2 - 1)! (j_1 + j_2  - 1)} \\
I_3 & \sim \frac{g(0) }{j_1! j_2!} (\log{M})^{j_1 + j_2}.
\end{align*}

In summary, we have shown so far that $I_1 \sim c (\log{X})^{j_1 + j_2}$ for some positive constant $c$ (the positivity follows by inspection of \eqref{eq:I2} and \eqref{eq:I3}).

We now return to the estimation of $I_\gamma$.  We move the line of integration over $v_1$ to a contour of the type $\gamma_r$ but which is twice as close to the imaginary axis as $\gamma_2$, thereby picking up a pole at $v_1 = -v_2$ only.  The integration over the new contour is small (arbitrary power of $\log{X}$ savings) by the same argument that showed $I_{1, \gamma_1}$ is small.  

The contribution from the pole at $v_1 = -v_2$ is
\begin{multline}
\label{eq:weirdpole}
-\frac{1}{(2 \pi i)^3} \int_{\gamma_2} \int \int \frac{\zeta(1 + 2s_1)\zeta(1 + 2s_2)\zeta(1 + s_1 + s_2)}{\zeta(1 + s_1 - v_2)\zeta(1 + s_1 + v_2)\zeta(1 + s_2 - v_2)\zeta(1 + s_2 + v_2)} 
\\
\cdot
\frac{M_1^{s_1}}{s_1^{j_1 + 1}}
\frac{M_2^{s_2}}{s_2^{j_2 + 1}} 
\frac{(V_2/V_1)^{v_2}}{v_2^2}
%\eta(s_1, s_2, -v_2, v_2)
 %G(-v_2) G(v_2) \Gamma(-v_2 + 1) \Gamma(v_2 + 1) ds_1 ds_2 dv_2.
g(s_1, s_2, -v_2, v_2) ds_1 ds_2 dv_2.
\end{multline}

We get different behavior in this integral depending on whether $V_2 > V_1$ or $V_2 = V_1$.  First assume $V_2 > V_1$.  In this case
we simply move $s_1$ and $s_2$ to contours very close to the imaginary axis as in the proof of the bound on $I_{1,\gamma_1}$.

Now assume $V_1 = V_2 = V$ and $M_1 = M_2 = M$ and continue with \eqref{eq:weirdpole}.  As a shortcut we may use Proposition \ref{prop:integralone} to compute the integrals over $s_1$ and $s_2$ because the $\zeta^{-1}(1 + s_i \pm v_2)$ factors do not have any poles in the regions of integration.  
We have that \eqref{eq:weirdpole} asymptotically evaluates as 
\begin{equation*}
- \frac{1}{4} \frac{(j_1 + 1)! (j_2 + 1)!}{j_1 + j_2 + 3} (\log{M})^{j_1 + j_2 + 3} \frac{1}{2 \pi i}\int_{\gamma_2} \frac{1}{\zeta^2(1 + v_2)\zeta^2(1  - v_2)} 
\frac{1}{v_2^2}
%\\
%\cdot
g(0, 0, -v_2, v_2) dv_2.
\end{equation*}
The integrand is holomorphic at $v_2 = 0$ so we may move the line of integration to the line $\text{Re } v_2 = 0$, for which it is obvious that the above expression is positive.
\end{proof}

\end{document}